%
%
%
%

\input harvmac.tex
\hfuzz 15pt



\def\za{\alpha} \def\zb{\beta} \def\zg{\gamma} \def\zd{\delta}
\def\ze{\varepsilon}   
\def\zk{\kappa} \def\zl{\lambda} \def\zm{\mu} \def\zn{\nu}
  \def\zr{\rho} \def\zs{\sigma} 
\def\zt{\tau}

\def\bL{\bar{\Lambda}}
\def\zG{\Gamma}   
\def\zL{\Lambda}  

\def\IZ{Z\!\!\!Z}
\def\[{\,[\!\!\![\,} \def\]{\,]\!\!\!]\,}
\def\dC{C\kern-6.5pt I}

\def\bw{\bar w}
\def\by{\bar y}

\def\bW{\overline W}
        \def\CC{{\cal C}}
        \def\CF{{\cal F}}
\def\CG{{\cal G}}

\def\CP{{\cal P}}        \def\CR{{\cal R}}

\def\un{{\bf 1}}

%

\def\({ \left( }
\def\){ \right) }

\catcode`\@=11
\def\Eqalign#1{\null\,\vcenter{\openup\jot\m@th\ialign{
\strut\hfil$\displaystyle{##}$&$\displaystyle{{}##}$\hfil
&&\qquad\strut\hfil$\displaystyle{##}$&$\displaystyle{{}##}$
\hfil\crcr#1\crcr}}\,}   \catcode`\@=12



\def\IC{\relax\hbox{$\inbar\kern-.3em{\rm C}$}}

\input amssym.def
\input amssym.tex
\def\IZ{\Bbb Z}\def\IC{\Bbb C}
\def\gg{\goth g} 


%
%

\def\go{{\goth g}} \def\bgo{\overline{\goth g}} 
  \def\bh{\bar{\goth h}}      
\def\th{\tilde{\goth h}}

\def\bR{\overline\Delta}                   
\def\Rr{\Delta^{\rm re}} \def\bP{\overline\Pi}
\def\br{\overline\zr}                      

\def\WL{P} \def\RL{Q}                      
 
\def\FW{\Lambda}\def\fw{\overline\FW}       

\def\bW{\overline W} \def\tW{\tilde W}      
\def\hw{w} \def\bw{\overline w}

\def\bKW#1{\bW^{[#1]}}                       
\def\unit{e}                                 
\def\UNIT{\blacktriangleleft\kern-.6em\blacktriangleright}

              
 \def\tr{t} \def\zt{\tau}
\def\ts#1,#2{{\tt e}^{#1\zk#2}}

\def\AA{A}                
\def\aa{a} 
 \def\gg{\gamma}

\def\iii{\iota}                           

\def\CC{{\cal W}^{(+)}} \def\tC{\tilde\CC} 

\def\UU{{\cal U}}                         

\def\supp{{\rm supp}}                          
\def\HT {{\cal T}}         
\def\VV{{\cal V}}   
\def\GG{{\cal G}}                 
\def\VM{V}                          

\def\cc{\chi} \def\bc{\overline\cc}       
\def\ch{{\rm ch}}

\def\dd{d_{\zk}}                           
\def\KK { K}                                 
\def\bK{\overline{
K}}
\def\ml{m}                                 
\def\bm{\overline{m}}                 
\def\ff{f}                               
\def\FF{F}

\def\CR{{\goth W}} \def\tCR{\tilde{\goth W}}
\def\gV{\goth V}                    

\def\fCC#1{\CC_{\le #1}} \def\gCC#1{\CC_{= #1}} 

\def\IZp{\IZ_{\ge 0}} 
\def\IZpp{\IZ_{>0}}
\def\IQ{{\Bbb Q}}




\def\la{\langle} \def\ra{\rangle}
\def\bs{\backslash} 

\def\ol#1{\overline{#1}}

\def\PROP{\medskip\noindent{\bf Proposition}\ \ }
\def\REMARK{\medskip\noindent{\it Remark.}\quad}
\def\LEMMA{\medskip\noindent{\bf Lemma}\ \ }
\def\PROOF{\smallskip\noindent{\it Proof:} \quad}
\def\endPROOF{\quad$\square$\medskip}
\def\COR{\medskip\noindent{\bf Corollary}\ \ }



\nref\AY{Awata, H., and  Yamada, Y.: 
       Fusion rules for the fractional level $\widehat{sl}(2)$ 
algebra.  
        Mod. Phys. Lett. {\bf A7}, 1185 - 1195 (1992)}

\nref\BI{  Bannai, E., and  Ito, T.: {\it Algebraic Combinatorics I: 
       Association Schemes .} Benjamin/Cummings, 1984}

\nref\Bour{ Bourbaki, N.: {\it Groupes et alg\`ebres de Lie.}
           Paris: Hermann, 1968}

\nref\BG{ De Boer, J. and  Goeree, J.: 
    Markov traces and II(1) factors in conformal field theory.  
    Comm. Math. Phys. {\bf 139},  267 - 304 (1991)}

\nref\DFZ{Di Francesco, P., and Zuber, J.-B.: 
        SU(N) lattice integrable models associated with graphs.
          Nucl. Phys. {\bf B338}, 602 - 646 (1990) 
  \semi
Di Francesco, P., and Zuber, J.-B.: SU(N) lattice integrable models 
and 
modular invariance,
in {\it Recent Developments in Conformal Field Theories}, Trieste
Conference 1989, Randjbar-Daemi, S., Sezgin, E., and Zuber,
J.-B., eds. (World Scientific 1990) 
\semi
Di Francesco, P.:
          Integrable lattice models, graphs and modular 
            invariant conformal field theories.
          Int. J. Mod. Phys. {\bf A7}, 407 - 500 (1992)}

\nref\FMa{ Feigin, B.L., and Malikov, F.G.: 
        Fusion algebra at a rational level and cohomology of 
           nilpotent subalgebras of $\widehat{sl}(2)$.
      Lett. Math. Phys. {\bf 31}, 315 - 325 (1994)} 

\nref\FM{ Feigin, B.L., and Malikov, F.G.: 
     Modular functor and representation theory of 
      $\widehat{sl}(2)$ at a rational level, in
         {\it Operads: Proceedings of Renaissance Conferences},  
     Cont. Math. {\bf 202}, p. 357, Loday, J.-L., Stasheff, J.D., 
	 and 
     Voronov, A.A., eds. (AMS, Providence, Rhode Island 1997)}

\nref\FGPa{ Furlan, P., Ganchev, A.Ch., and Petkova, V.B.: 
             Quantum groups and fusion rule multiplicities.
   Nucl. Phys. {\bf B343}, 205 - 227 (1990)} 
\nref\FGP{ Furlan, P., Ganchev, A.Ch., and Petkova, V.B.:
        Fusion rules for admissible representations of affine 
algebras:   
        the case of $A_2^{(1)}$. 
         Nucl. Phys. {\bf B518} [PM], 645 - 668 (1998)}
                                              
\nref\Hump{ Humphreys, J.M.: {\it Reflection Groups and
           Coxeter Groups}, Cambridge: Cambridge University Press, 
1990} 

\nref\K{ Kac, V.G.: {\it Infinite-dimensional Lie Algebras}, third
        edition, Cambridge: Cambridge University Press, 1990}

\nref\KW{ Kac, V.G., and  Wakimoto, M.:  
      Modular invariant representations of infinite-dimensional 
      Lie algebras and superalgebras. 
    Proc. Natl. Sci. USA {\bf 85}, 4956 - 4960 (1988)   \semi 
 Kac, V.G., and  Wakimoto, M.:
         Classification of modular invariant 
                 representations of affine algebras.
               Adv. Ser. Math. Phys. vol {\bf 7}, pp. 138 - 177. 
               Singapore: World Scientific 1989
         \semi
 Kac, V.G., and  Wakimoto, M.: 
        Branching functions for winding subalgebras and tensor 
products.
          Acta Applicandae Math. {\bf 21}, 3 - 39 (1990)}

\nref\Pas{ Pasquier, V.: Operator content of the ADE lattice models.
 J. Phys. {\bf A20}, 5707 - 5717 (1987) }
   
\nref\PZ{ Petkova, V.B.,  and Zuber, J.-B.: 
          From CFT to graphs.
          Nucl. Phys. {\bf B463}, 161 - 193 (1996) }

\nref\V{ Verlinde, E.:
      Fusion rules and modular transformations in 2D conformal field 
theory.
           Nucl. Phys. {\bf B300} [FS22], 360 - 376 (1988) }

\nref\MW{ Walton, M.: 
       Fusion rules in Wess-Zumino-Witten models.
           Nucl. Phys. {\bf B340}, 777 - 790 (1990) }


\null
\vskip 3.5cm

\centerline{\bf  An Extension of the Character Ring of sl(3)
and Its Quantisation}\footnote{}{furlan@ts.infn.it, 
ganchev@inrne.bas.bg, petkova@ts.infn.it}

\vskip 1cm
\centerline{P. Furlan$^{*, \dagger}\, $, $\ $  A.Ch.
Ganchev$^{** }$ $ \ $
 and $\ $ V.B. Petkova$^{**, \dagger}$} 
\bigskip
\bigskip
\centerline{\it $^{*}$Dipartimento  di Fisica Teorica
 dell'Universit\`{a} di Trieste}
\centerline{\it Strada Costiera 11, 34100 Trieste, 
 Italy}
\centerline{\it $^{\dagger}$Istituto Nazionale di Fisica Nucleare
(INFN), Sezione di Trieste, Italy} 

\medskip

\centerline{\it $^{**}$Institute for Nuclear Research and
Nuclear Energy} 
\centerline{\it Tzarigradsko Chaussee 72, 1784 Sofia, 
Bulgaria}
\medskip

\bigskip\bigskip

\vskip 3cm

\centerline{\bf Abstract}
\medskip
\medskip
\medskip

We construct a commutative ring with identity which extends the
ring of characters of finite dimensional representations of
$sl(3)$. It is generated by characters with values in the group
ring $\IZ[\tW]$ of the extended affine Weyl group of
$\widehat{sl}(3)_k$ at $k\not \in \IQ$.  The `quantised' version at
rational level $k+3=3/p$ realises the fusion rules of a WZW
conformal field theory based on admissible representations of
$\widehat{sl}(3)_k\,.$

\vfil\eject

%
%
\newsec{Introduction}

The aim of this work is to describe the characters
(one-dimensional representations) of the fusion algebra  of 
$\go=\widehat{sl}(n)_k$ WZW conformal models at rational 
(shifted) level $\kappa:=k+n=p'/p\,,$ for $n=3=p'\,, $ 
 $(p,3)=1$ and their
`classical' counterparts.

In \FGP\ we have realised this fusion algebra as a matrix algebra
$\CF\subset$ Mat$_{p^2}({\cal C})$ with integer nonnegative
structure constants ${}^{(p)}N_{y, x}^z\,$ and a basis $\{N_y\,,
\ (N_y)_x^z={}^{(p)}N_{y, x}^z\,, \ \ N_{\un}=\un_{p^2}\}\,.$  The
labels of the basis run over the highest weights of the
admissible representations at $\kappa=3/p$ \KW, which
conveniently are parametrised  by a subset $\CC_p$ of the affine
Weyl group $W$,
see the text for precise definitions.  $\CF$ is a
matrix realisation of a commutative, associative algebra with
identity, a distinguished basis, and an involution $*$.  The
definition of the fusion algebra -- an example of $C$-algebra
(Character - algebra)  in the terminology of \BI{}, implies that
$N_y$ are normal, hence simultaneously diagonalisable by a
unitary matrix $\psi_y^{(\mu)}$ labelled by some set
$E_p=\{\mu\}$ of $p^2$ indices, 
\eqn\pasv{
{}^{(p)}N_{x,y }^z=
\sum_{\mu\in E_p}\, {\psi_{x}^{(\mu)}\over\psi_{\un}^{(\mu)}}\
\psi_y^{(\mu)}\, \psi_z^{(\mu)\,*}
=\sum_{\mu\in E_p}\,
\cc_x^{(p)}(\mu)\ {\cc_y^{(p)}(\mu)\ \cc_{z^*}^{(p)}(\mu)\over 
\sum_{u\in \CC_p}\, |\cc_u^{(p)}(\mu)|^2 }\,.
}
The eigenvalues
$\cc_{x}^{(p)}(\mu)=\psi_x^{(\mu)}/\psi_{\un}^{(\mu)}$ of  $N_x$
provide $p^2$  linear representations $N_y\rightarrow
\cc_y^{(p)}(\mu)\,, $ i.e., characters of $\CF$, labelled by the
set $E_p$,
\eqn\Iqch{
  \cc_{x}^{(p)}(\mu)\ \cc_{y}^{(p)}(\mu)
= \sum_{z\in \CC_p} {}^{(p)}N^{z}_{x,y} \ \cc_{z}^{(p)}(\mu)\,. 
}

With an appropriate 
reinterpretation of 
the labels this is the general setting for any RCFT.
In particular for WZW models based on the $\widehat{sl}(n)_k$ 
integrable representations (a subclass of the 
admissible representations at integer level $k=p'-n$),
both indices of  the unitary matrix $\psi_{\zl}^{(\mu)}$ belong
to the alcove 
$P_+^{k}\,,$ $\, \psi_{\zl}^{(\mu)}$ 
is a symmetric matrix which coincides with the modular matrix, 
$\psi_{\zl}^{(\mu)}=
S_{\zl \zm}^{(p')}$, while  \pasv\ reduces to the Verlinde
formula \V\ for the fusion rule (FR) multiplicities, 
i.e., the dimensions of the spaces of chiral vertex operators.

In \FGP\ we have specified the above general setting to the case
of the generic subseries $\kappa=3/p\,$ of the admissible
representations, describing explicitly one of the fusion
matrices, a `fundamental' fusion matrix $N_f$. Here $f$ is some
analog of the $\bgo= sl(3)$ fundamental weight $\bL_1=(1,0)$, but
the  product of $N_f$ with any  $N_x$ produces generically seven
terms, each appearing with multiplicity one.  The numbers
${}^{(p)}N_{f,y}^z=0,1$ were found by solving a set of  algebraic
equations coming from the decoupling of singular vectors in
$\go=\widehat{sl}(3)_k$ Verma modules. More precisely the
decoupling of the `horizontal' singular vectors (which exist also
for generic values of the level, $\kappa\not\in \IQ$) determines
the generic seven points, while the additional truncation
conditions, including the ones at rational level, were obtained
making also some assumptions suggested from explicit computations
at small values of $p$.  Although  $N_f$ together with its
conjugate $N_{f^*}$ are  not sufficient to build up the full
polynomial fusion ring, diagonalising this matrix  in an
orthonormal basis, i.e., finding the  eigenvector matrix
$\psi_y^{(\mu)} $  (common to all fusion matrices) recovers
through \pasv\ all  FR multiplicities.  We called  \pasv\
Pasquier--Verlinde type formula since similar formul\ae\ --
though with different interpretation of the structure constants
(no more required to be nonnegative integers), were  first
discussed in \Pas\ in the context of lattice ADE models. Such
formul\ae, together with their dual counterparts, describing the
structure constants of a `dual' algebra, have been furthermore
exploited in the study of nondiagonal modular invariants of the
integrable WZW conformal models \DFZ, \PZ.

Diagonalising a $p^2 \times p^2$ matrix becomes a tedious task
for big $p$ (we have done this exercise for $p=4,5$) so it is
preferable to have an explicit analytic formula for $\psi_y^{(\mu)}
$, or equivalently, for the  characters $\cc_y^{(p)}(\mu)$ of the
fusion algebra.  The path we follow in this paper to find the
characters $\cc_y^{(p)} (\mu)$ was suggested by the second
formula for the admissible FR multiplicities  conjectured (with
slightly changed notation) in \FGP,
\eqn\ws{
{}^{(p)}N^{z}_{x,y}: = \sum_{w'\in W^{[z]}}
\,\det(\overline{w'})\ \ml^{x}_{w'\,z\, y^{-1}}\,.  
}
This formula  is analogous to the formula for the integrable FR
multiplicities derived in  \K,\MW, \FGPa, which is a truncated
version of the classical Weyl--Steinberg formula for the tensor
product multiplicities of finite dimensional representations of
the horizontal subalgebra $\bgo$ (with the role of the horizontal
Weyl group $\bW$ taken by the affine Weyl group $W$). In \ws\ the
summation runs over the Kac--Wakimoto (KW) affine group
$W^{[z]}$, generated by   Kac--Kazhdan reflections corresponding
to  singular vectors of the $\go$ Verma module of highest weight
parametrised by the element $z\in W$. In the integrable case the
counterpart of $m^x_v$ is the multiplicity $m^{\zl}_{\zm}$   of
the weight $\mu$ of the $\bgo$ finite dimensional module of
highest weight $\zl$.
Here the integers $ \ m^x_v\,, x,v \in W\,$ 
describe a
generalised (finite) weight diagram.  The realisation that one
has  to generalise the classical notion of weight diagram in
order to describe the rational level FR was the main lesson of
the, otherwise still incomplete,  analysis of the
null-decoupling equations of \FGP. While the $\bgo$ 
weight diagrams 
are  subsets of the root lattice $Q$ of  $\bgo$, attached to the
highest weight $\zl$, in our case the role of a `root lattice' is
taken by the affine Weyl group $W=\bW \ltimes t_Q \,$ at a generic
($\kappa \not \in \IQ$) level.  Now the question is, can one find
 formal  characters encoding this information about
the generalised weight diagrams, which furthermore are closed
under multiplication and recover  a `classical' analogue of
\ws\ at generic levels, with the affine KW groups replaced by
their horizontal counterparts $\bW^{[z]}$.  Then by analogy with
the integrable case one can `quantise' these `classical' characters
imposing the periodicity conditions, accounting for the rational
level, and thus recover the linear representations of $\CF$.

Though it is not necessary, to understand our way of reasoning
and the motivating idea it is helpful to imagine that there is a
finite dimensional algebraic object playing the role of $\bgo$.
The irreps of the assumed hidden  algebra are labelled  
in general by
$x\in\tC$, a fundamental subset of the extended affine Weyl group
$\tW=\bW\ltimes t_P \,$ with respect to the right action of $\bW$
(equivalent to the action of the horizontal KW Weyl group, see
\KW\ and the text below).  The same set  also parametrises
maximally reducible Verma modules of $\go$ at generic level and
the corresponding irreducible quotients.  In the simplest
$\widehat{sl}(2)_k$  case solved in \AY, \FMa\ this intrinsic algebra
is the $\IZ_2$ graded algebra $osp(1|2)$ \FM; 
the Weyl group $\bW\simeq \IZ_2$ of $sl(2)$ distinguishes between 
even and odd vacuum state finite dimensional
irreps of $osp(1|2)$.   Alternatively
mapping $\iii:\tW \rightarrow Q$ one can use the horizontal
subalgebra $\bgo=sl(2)$ itself, but restricting to its
representations with highest weights on the root lattice $Q$
(integer isospins).  In general this map $\iii$ (see section 2)
allows to describe the generalised modules through subsets of the
supports of standard  modules of $\bgo$ (with  highest weights of
$n$-ality zero); the case $sl(2)$ is trivial since then
Im$(\iii) \equiv Q$.

The construction of  the  `classical characters' at generic
level, $k\not \in \IQ$, takes a considerable part  of this work
and is its main novel result. 
It is done in the first part of the paper
(sections 2,3,4) by full analogy with the standard case starting
with 
some  analogs of the Verma module characters, 
used as ingredients in a generalised Weyl character formula.
While the input about the relevant
generalised supports of the finite or infinite `modules' (the set
of weights and their multiplicities) is essentially taken over
from \FGP, slightly rephrased and generalised, the main effort
here is to find the proper multiplicative structure of the formal
characters.  They are  elements of the group ring $\IZ[\tW]$ of
the extended affine Weyl group  of $\widehat{sl}(3)_k$, closed
under multiplication. The structure constants of the resulting
commutative subring $\tCR$ of $\IZ[\tW]$ satisfy  a generalised
Weyl-Steinberg formula.  As a side result one obtains an explicit
formula for the cardinality of the generalised weight diagrams,
i.e., for  the dimensions of the `finite dimensional  modules' of
the unknown algebra possibly generalising $osp(1|2)$.  
Some steps of the construction in this first part hold, or are
straightforwardly generalisable, for $\go=\widehat{sl}(n)_k\,,$
arbitrary $n$. So we keep the exposition general although 
we present the details fully for the case
$\go=\widehat{sl}(3)_k$; the general case will be elaborated
elsewhere.  In the second part (sections 5,6,7) we consider
rational values $\kappa=3/p$, thus a generic subseries of the
admissible representations of $\go=\widehat{sl}(3)_k$. The formal
characters are `quantised' imposing periodicity constraints and
realised as $\IC$ - valued functions.  They satisfy the
orthogonality and completeness relations equivalent to the
unitarity of $\psi_y^{(\mu)}$, so that one  recovers the
Pasquier--Verlinde type formula \pasv.  Thus given the characters
both  formul\ae\ for the FR multiplicities of the admissible
representations \pasv\ and \ws\ derive from the classical analog
of the Weyl--Steinberg formula established in the first part of
the paper. Furthermore one proves a third formula for the FR
multiplicities, conjectured in \FGP, which also has a `classical'
counterpart,
\eqn\qstrci{
{}^{(p)}N^{z}_{x,y} ={}^{(3 p)}\bar{N}^{\iii(z)}_{\iii(x)\,
\iii(y)} \,, 
}
where ${}^{(3 p)}\bar{N}^{\iii(z)}_{\iii(x)\,\iii(y)}$ are
structure constants of the integrable  fusion algebra at
$\kappa=3p$, while $\iii$ is the map sending $\tW$ to a subset of
the triality zero weights. Presumably the fusion rules of
$\widehat{sl}(n)$ WZW at  $\kappa =n/p$ are given by the above
formula with $n$ substituting $3$.

%
%
\newsec{Preliminaries.}

We start with fixing some notation.  Let  $\bR\,,$ $\bR_+\,,$
$\bP=\{\za_1,\dots, \za_{n-1}\}\,$  be, respectively, the sets of
roots, positive roots, and the simple roots of  $\bgo=sl(n)$, and
$\Rr\,,$ $\Rr_+=\bR_+\ \cup\ (\bR+\IZpp\ \zd)\,,$
$\Pi=\{\za_0=\delta- \sum_i\za_i\,,\, \bP\}\,,$ -- the set of
real roots, real positive roots and the simple roots  of the
affine algebra $\go=\widehat{sl}(n)_k$.  Let $\fw_i$ be the
fundamental weights of $\bgo$, i.e.,
$\la\fw_i\,,\za_j\ra=\zd_{ij}$, with respect to the
Killing-Cartan bilinear form $\la\cdot\,,\cdot\ra$ on the dual
$\bh^*$ of the Cartan algebra.  With $\RL=\oplus_{i}\ \IZ\ \za_i$
and $\WL=\oplus_{i}\ \IZ\ \fw_i$ we denote the root and weight
lattices of $\bgo$. Their positive cones are $\RL^+=\oplus_{i}\
\IZp\ \za_i$ and $\WL_+=\oplus_{i}\ \IZp\ \fw_i$.  The negative
cone $-\RL^+$ is the support of the $\bgo$ Verma module of 0
highest weight, while $\WL_+$ is the chamber of integral dominant
weights -- the highest weights of the finite dimensional
representations of $\bgo$.  The form $\la\cdot\,,\cdot\ra$
extends to $\th^*=\bh^*\oplus\IC\FW_0\oplus \IC\zd$ with $\la
\bh^*, \IC\FW_0\oplus \IC\zd \ra=0\,,$ $\la \FW_0\,,
\FW_0\ra=0=\la \zd\,,
\zd\ra\,,$ $\la\zd\,,\FW_0\ra=1$.  The fundamental weights of
$\go$ are $\{\FW_i=\FW_0+\fw_i\}\,,$ $\fw_0=0\,.$ The
`horizontal' projection $\overline{(\cdot)}:\th^*\to\bh^*$ is
defined as having the kernel $\IC\zd\oplus\IC\FW_0$.  The Weyl
vector is $\zr=\sum_j\, \FW_j\,.$ The   Weyl group  $\bW$ of
$\bgo$ is a finite Coxeter group generated by the simple
reflections $\hw_i=\hw_{\za_i}$ with relations
$(\hw_i)^2=1=(\hw_i\hw_{i+1})^3$, $\hw_i\hw_{j}=\hw_j\hw_{i}\,,$
$j\not= i\pm 1\,, \ i,j=1,2,\dots, n-1$. The affine Weyl group
$W$ is 
generated by the simple reflections $w_j\,, \ j=0,1,\dots, n-1$
with the same type of relations (for $n>2$), identifying
$w_n=w_0$.  These groups can be depicted by their Cayley graphs,
in which the vertices correspond to elements of the group, edges
to the generators, and the elementary polygons to the relations.
E.g., fig. 1 depicts (a finite part of) the  Cayley graph of the
affine Weyl group $W$ of $\go=\widehat{sl}(3)$, presented as $\{
\hw_i :  (\hw_i)^2=1=(\hw_i\hw_j)^3 , \{i,j\}\subset\{0,1,2\}
\}\,, $ while  any of the `12' elementary hexagons on fig. 1 is
the Cayley graph of the corresponding horizontal Weyl group
$\bW$.  The labels $i$ on the edges correspond to the generators
$\hw_i$, $i=0,1,2$. The `origin' of the graph, the vertex
corresponding to the group unit $\un$ is denoted by $\UNIT$. The
three types of hexagons (`12', `01', and `20') correspond to the
three Artin type relations among the generators.  It is
convenient to introduce the following shorthand notation:
$\hw_{ijk\dots} = \hw_i \hw_j \hw_k \dots$.

We will also need the extended affine group $\tW$ 
defined as the semi-direct product
$\tW=\bW\ltimes t_P\,$ (while $W=\bW\ltimes t_Q\,$),  $t_P$ being
the subgroup of translations in the weight lattice $P$. The
elements $t_{\beta}\,,$ $\beta\in P\,,$ act on $\th^*$ as
\eqn\Iee{
\tr_{\zb}(\zL)
 = \zL + \la\zL\,,\zd\ra\zb - \Big(\la\zL\,,\zb\ra + 
   {1\over 2}\la\zb\,,\zb\ra\la\zL\,,\zd\ra\Big) \zd \,,
}
and for $ \za\in \bR\,,\,\, l\delta-\za\in \Rr_+\,, \,\, y\in
\tW\,,\,\, \zb\in P\,,$ one has the properties
\eqn\Ie{
 t_{l\za}=w_{l\delta-\za}\,w_{\za}\,, \quad 
   y\,t_{\zb}\,y^{-1} = t_{\bar{y}(\beta)}\,, \quad
    y\, w_{\za}\,y^{-1}=w_{y(\za)}\,. 
}
We have denoted by overbar the projection of $\tW$ onto the
horizontal subgroup $\bW$ sending the affine translations to the
unit element.

The group $\tW$ can be also written as $\tW=W\rtimes A\,,$ where
$A$ is the subgroup of $\tW$ which keeps invariant the set of
simple roots $\Pi$ of $\go$. It is a cyclic group generated by
$\zg= t_{\fw_1}\, \bar{\zg}$, where $ \bar{\zg}=\bw_1\dots
\bw_{n-1}$ is a Coxeter element in $\bW$ generating  the 
cyclic subgroup $\bar{A}$ of $\bW$.  One has
$\zg(\za_j)=\za_{j+1}=\zg^{j+1}(\za_0)$ for $j=0,1,2,\cdots, n-1$
identifying $\za_n\equiv\za_0$.  In the case of $\bgo=sl(3)$ we
will think of the Cayley graph of $\tW$ as a 3-sheeted covering
of the graph of $W$ with, for example, the ``fiber'' over the
edge `0' connecting the vertices $\un$ and $\hw_0$ being the set
$\UU=\{A, A\,\hw_0\}$ and this part of the graph of $\tW$ is
depicted on fig. 2. The oriented edges correspond to $\zg$ and
the squares -- to the implementation of the automorphism  of
$W\,,$ $w_{\za}\rightarrow \zg \, w_{\za}
\zg^{-1}=w_{\zg(\za)}\,,$ $\za\in \Pi$.

Introduce the set $\CP=\{\zL= y\cdot k \zL_0\,, \,\, y\in \tW \}$
where the  shifted action of 
$\hw\in \tW$ 
on $\th^*$ is given by 
$\hw\cdot\zL=\hw(\zL+\zr)-\zr$.  In this and the following two
sections we shall assume that the level $k$ is generic, $k \not
\in \IQ$, which in particular ensures that if $y\cdot k \zL_0=k
\zL_0$ then $y\equiv \un$.

According to the general criterion of Kac-Kazhdan if
$\la\zL+\zr,\zb\ra $ is a positive integer for some $\zb\in
\Rr_+$ the Verma module $M_{\zL}\,,$  $\zL\in \CP$, of  $\go$ is
reducible, containing a Verma submodule $M_{w_{\zb}\cdot \zL}\,.$
In particular if for some $y=\by\, t_{-\zl}\in \tW\,$ and $
\za\in \bR_+\, $ we have $y(\za)=\la\zl,\za \ra\, \zd+\by(\za)\in
\Rr_+\,,\ $ then the KK condition is fulfilled for  $\zL=y\cdot
k\zL_0$ and the root  $\zb=y(\za)$, and using \Ie\
\eqn\Ia{
w_{y(\za)}\cdot \zL=y\,w_{\za}\,y^{-1}\cdot \zL=y\, w_{\za}\cdot
k\zL_0 = \by\,w_{\za} t_{-w_{\za}^{-1}(\zl)}\cdot  k\zL_0 \,.  
}
If $y(\za_i)\in \Rr_+\,,$ $\forall$ $\za_i\in \bar{\Pi}\,,$ the
reflections $\hw_{y(\za_i)}$ generate an  isomorphic to $\bW$
group $\bKW{\zL}$ (to be denoted also $\bKW{y}$), introduced by
Kac-Wakimoto \KW. We shall refer to these groups as  (horizontal)
KW (Weyl) groups.  
According to  \Ia\ we can identify the action of a KW group with
the right action of $\bW$ on $\tW$.

Now we introduce the subset $\CP_+\subset\CP$ of weights $\zL$
such that the corresponding Verma modules $M_\zL$ are maximally
reducible. {}From the above discussion it is clear that $\
\CP_+=\tC\cdot k\zL_0\,,$ where
\eqn\Ida{
\tC=\{y\in \tW\,|\, y(\za_i)\in \Rr_+\ \ {\rm for}\ \forall
\za_i\in \bP\}\,. 
}
 The Bruhat ordering on the orbit
$\bKW{\zL}\cdot k\zL_0$ describes the embedding pattern among the
Verma modules $\{ M_{\zL'},
\zL'\in  \bKW{\zL}\cdot k\zL_0 \}$. 

Denote $\CC=\tC\cap W$. One has
\PROP {\bf 2.1\ } {\it $\CP_+$ is a fundamental domain in $\CP$  with
respect to the action of the KW Weyl groups, or, equivalently, $\
\tC$ ($\CC $) is a fundamental domain in $\tW $ ($W$) with
respect to the right action of $\bW \,.$}
\PROOF Let us introduce some notation first. Let   $H_{\za}=
\{\zl \in \bh^*\,| \, \la \zl\,,\, \za\ra=0  \}$ be
the hyperplane orthogonal to the root $\za$. 
 For $w\in\bW$ let $I(w)=\{i \,|\, w(\za_i)<0\,,\, \za_i\in\bP\}$
and  
$$
P_+^{(w)}:= \{ \zl \in P_+ |\, \la\zl,\za_i\ra>0\,,\, 
i\in I(w) \} = P_+\bs \Big(\big(\cup_{i\in I(w)}\,
H_{\za_i}\big)\cap P_+\Big) = P_+ +\sum_{i\in I(w)}\,\bL_i \,.
$$
The definition \Ida\ can be obviously rewritten as 
\eqn\Id{
\tC =\{y \in \tW\,|\, y=w t_{-\zl}\,,
\,\,  w\in \bW\,, \,\, \zl \in P_+^{(w)} \} \,.
}
Hence the statement of the Proposition follows from $
\cup_{_{w\in \bW }} \tC\, w= \cup_{_{w\in \bW }}\,\cup_{_{w'\in 
\bW }}\, t_{-w'(P_+^{(w')})}\,w'\, w = \cup_{_{w\in \bW }}\,
t_{-P}\, w = \tW\,$ once the following lemma is established:

\LEMMA {\bf 2.2\ }
{\it $P=\cup_{_{w\in \bW }} w(P_+^{(w)})$ is a partition, i.e., a
disjoint union.}

\noindent
To prove the lemma one has to exploit several standard properties
of the Weyl group $\bW$ which can be found e.g., in \Hump,
chapter I.
\PROOF 
Any $\zl' \in P$ is represented as $\zl'=w(\zl)$ for some $\zl\in
P_+$, $w\in \bW$. Denote $X=X_{\zl}=\{i\,|\, \la \zl \,, \za_i\ra=0
\}$.  According to Proposition 1.10c of \Hump\ the element $w$
splits uniquely into a product of two elements of $\bW$, $w=uv$,
s.t. $v\in W_X\,$ (the group generated by the simple reflections
labelled by the subset $X$), and $u(\za_i)>0\,,$ for any
$i\in X$. We have $v(\zl)=\zl\,,$ $\ I(u)\cap X = \emptyset\,,$ hence
$\zl\in P^{(u)}_+$ and $\zl'=u(\zl)\in u(P^{(u)}_+)$.  This
proves that $P$ is covered by the union of subsets
$w(P^{(w)}_+)$.  The uniqueness of the above splitting proves
also  the disjointness.
\endPROOF
\medskip

We shall refer to $\tC$ (or, equivalently, $\CP_+$) as a
`dominant chamber'. The left action of the group $A$ (the shifted
action of $A$) keeps  $\tC$ ($\CP_+$)  invariant and $\cup_{a\in
A} \, a\, \CC= \tC$. Similarly $\CP_+$ is $\IZ/n
\IZ$ graded by the $n$-ality $\tau(\zL=\by t_{-\zl}\cdot k \zL_0):=
\tau(\zl)\,,$ where $\tau(\zl)=\sum_i\, i \la \zl, \za_i\ra=  n\,
\la \zl, \bL_{n-1}\ra$ mod $n$ is the standard grading in $P$.
\medskip
 
In the case $sl(3)$ the chamber $\tC$ can be also expressed
as $\tC= \UU\ \tr_{- P_+} $ in terms of the subset $\UU=\{ A, A\,
w_0\}\subset \tW\,,$ depicted on fig. 2. 

\medskip

Next we introduce a map $ \iii$ of $\tW$ (or $\CP$)  into $Q$ 
\eqn\Ib{
 \iii: \tW\ni y=\by t_{-\zl} \
                                \mapsto \
     n\,\zl + \by^{-1}\cdot 0=n\,\zl \ -
         \!\!\sum_{_{
          \zb>0\,,\ \by(\zb)<0 }}\!\!\beta \,\in Q\,.
}
See \K\  (exercise 3.12 of  ch. 3) for the last equality.
$A$ is mapped by $\iii$ to zero.

 The map $\iii$ has the ``twisted log'' property
\eqn\tlg{
\iii(x y ) = \ol{y}^{-1}(\iii(x)) + \iii(y) \,.
}
Compare with the horizontal projection map 
\eqn\hor{\eqalign{
h:\tW\ni  y&=\ol{y}t_{-\zl} \mapsto h(y)
= \ol{y\cdot k\zL_0} = \ol{y}\cdot(-\zk\zl) 
\in \zk P + \bW\cdot 0\,,\cr 
h(x y) &=h(x) + \bar{x}(h(y))
\,.
}}

The map \Ib\ provides another equivalent definition of
the chamber $\tC$. Indeed comparing with \Id\ and
using that $1\le \la \br\,, \, \za
\ra\le n-1\,$ for any $ \za\in \bR_+$ one easily checks
\eqn\If{
\tC=\{y\in \tW\, |\, \iii(y) \in P_+\}\,.
}

The relation \tlg\   implies that $\iii$ intertwines between the
KW action  (equivalent according to \Ia\ to the right action of
$\bW$ on $\tW$) and the ordinary shifted action  of $\bW$, i.e.,
\eqn\Ic{
 \iii(y\bw)=
 \bw^{-1} \cdot\iii(y)\,, 
\quad \bw  \in\bW\,. 
}

Both subsets of $Q$, the image Im$(\iii)$ and its complement
are invariant  under the shifted action of the Weyl group $\bW$. 
The $\bgo$ Verma modules of highest weight $\iii(y)$ are
reducible iff the corresponding $\go$  Verma modules of highest
weight $\zL=y\cdot k\zL_0$ are reducible and the pattern of
embeddings of submodules in both cases is identical.

\medskip
On fig. 3 we have illustrated the map ${1\over 3}\iii$ for the
case $\go=sl(3)$. It maps the even (under the gradation
$\det(w)=\pm 1$) elements of $W$ to $\WL$ and the odd elements to
$\WL+\theta/3$ ($\theta= \alpha_1+\alpha_2$) in such a way that
the vertices of the Cayley graph (if we make it into a `rigid'
geometrical graph by fixing the length of each edge to be of
length $\sqrt{2}/3$ assuming as usual that the roots of $sl(3)$
are of length $\sqrt{2}$) geometrically `sit' at the same places
as the vertices of the two lattices $\WL$ and $\WL+\theta/3$.
(See the figure.) In other words refining by $3$ the weight
lattice, $\fw_i \rightarrow \fw_i/3$, (or equivalently,
rescaling $\kappa \rightarrow \kappa/3$) the points of the Cayley
graph can be identified with a subset of the triality zero
weights $\{\zl=\sum_i\, 3 n_i \fw_i/3\}\cup \{\zl=\sum_i\, (3
n_i+1) \fw_i/3 \}$ in the refined lattice, while the ``excluded''
triality zero points $\zl=\sum_i\, (3 n_i-1) \fw_i/3\in
P-\theta/3$, correspond to the centers of the elementary
hexagons.

\noindent
\REMARK
The above analysis properly extends to the larger than $\CP$
region 
$$
\{\zL= y\cdot(\zl'+ k \zL_0)\,, \ \,\, y\in \tW\,,\  \, \zl'
\in P \,; \ \ 
k \not \in \IQ \}\,.
$$
 It contains a subset with $y\in \tC\,, \, \zl'\in P_+\,,$
providing highest weights of `maximally reducible' $\go$ Verma
modules.  For our purposes it is sufficient to choose $\zl'=0$
thus restricting to weights $\zL$ parametrised by $y \in\tW$. 
%
%
\newsec{Characters.}

Let us start by recalling the  supports and characters  of
ordinary Verma and finite dimensional modules of $\bgo$.
The latter characters are elements of the group ring
$\IZ[\tr_\WL]$ of the group of translations by the weight
lattice.  Keeping with tradition, the translations $\tr_{-\zl}$
from the previous section will be written as formal exponentials
$\ts-,\zl\,;$ $\kappa =-1 $ recovers the standard notation, see,
e.g., \Bour\ for standard definitions.
The character of the Verma module $\VM_{\zl}$ of highest weight
$\zl$ is given by 
\eqn\cch{
  \ch(\VM_{\zl})
     = \sum_{\zm\in \zl-Q^+} \bK_{\zm}^{\zl} \ts-,{\zm} 
     =\ts-,{\zl}\, \sum_{\zb\in\RL^+} \KK_{\zb} \ts,{\zb} 
     = \ts-,{\zl} \  \prod_{\za\in\bR_+} \ (1-\ts,{\za})^{-1} 
     = \ts-,{(\zl+\br)}\  \dd  \,, 
}
where the multiplicity $\bK_{\zm}^{\zl}$  of an weight $\mu$ is
expressed via the Kostant partition function
$\bK_{\zm}^{\zl}:=\KK_{\zl-\zm}\,\in\IZp\,,$ while $\dd$ is a
$\bW$ invariant (up to a sign) quantity $w(\dd)=\det(w)\, \dd$
for $w(\ts,{\zb}):=\ts,{w(\zb)}$.  The support of a module is the
set of weights $\mu$ of nonzero multiplicity $\bK_{\zm}^{\zl}$,
thus $\supp\VM_{\zl}= \zl-\RL^+$.

The irreducible finite dimensional modules can be resolved in
terms of Verma modules, i.e., each Verma module
$\VM_{\bw\cdot\zl}$, with $\bw\in\bW$ and $\zl\in\WL_+$, contains
submodules of weight $\VM_{\bw'\cdot\zl}$ for all $\bw'>\bw$ in
the Bruhat ordering of $\bW$ (in a convention in which $\un$ is
the smallest element) and grading $\bW$ by the reduced length of
words the Verma module inclusions organize in a BGG
(Bernstein--Gelfand--Gelfand) resolution.  For the characters
$\bc_\zl\,,$ $\zl\in\WL_+$ of irreducible finite dimensional
modules the BGG resolution gives immediately  the Weyl character
formula
\eqn\sfch{\eqalign{
  \bc_\zl
 = \sum_{\hw\in\bW} \det(\hw)\ \ch(\VM_{\hw\cdot\zl})& = 
 \dd\ 
 \sum_{\hw\in\bW} \det(\hw)\  \ts-,{\hw(\zl+\br)}
  = \sum_{\zm\in P
  } \bm^{\zl}_{\zm}\ \ts-,{\zm} \,,\cr
\bm^{\zl}_{\zm}& = \sum_{w\in\bW} \det(w)\ 
K_{w\cdot\zl - \zm}\,, 
}}
and from $\bc_{\un}=1$ the factor $\dd$ is expressed as
$1/\dd=\sum_{\hw\in\bW} \det(\hw)\  \ts-,{\hw(\br)}$.  The
support (weight diagram) is $\zG_{\zl} =\{\zm\in\WL\,|\,
\bm^{\zl}_{\zm}\ne  0\}$.  {}From \sfch\ it follows that
\eqn\ginv{ 
w(\bc_\zl)=  \bc_\zl\,, \ 
\ \   \bm^{\zl}_{w(\zm)}=\bm^{\zl}_{\zm}\,,
}
 and,  extending the first line of \sfch\ to $\zl\in P$, 
\eqn\winv{
    \bc_{w\cdot\zl}=\det(w)\,  \bc_\zl\,.
} 
The map $\iii$ introduced in \Ib\ establishes a correspondence
between (highest weights of) Verma modules of the affine algebra
$\go$ and of the horizontal subalgebra $\bgo$, with identical
reducibility structure.  Now we shall introduce another class of
`Verma modules' $\VV_y$ 
(of yet unknown finite dimensional algebra)
described through its supports, i.e., weights and their
multiplicities. The supports are parametrised by elements of $a W
\subset\tW$ (for a fixed $ a\in A$) and mapped by $\iii$ into
subsets of supports of  $\bgo$ Verma modules.  Motivated by  the
intertwining property \Ic\ we shall call such a module reducible
if its  $\iii$ image is a reducible Verma module. In particular
the maximally reducible `Verma modules' have highest weights
$\zL=y\cdot k\zL_0\,,$ $y\in \tC\,,$ whence the name `dominant
chamber' for the latter subset of $\tW$.  Furthermore in full
analogy with the representation theory of $\bgo$ we shall
introduce `finite dimensional modules' obtained factorising
maximal  `Verma submodules' of reducible `Verma modules' by an
analog of the BGG resolution in which the role of $\bW$ is
replaced by the action of the KW group, equivalent according to
\Ia\ to the right action of $\bW$ on $\tW$. We have already
explained on the example of $\bgo=sl(3)$ this construction in
\FGP\ on the level of multiplicities of weights, here we add a
realisation of the characters of these `Verma- and finite
dimensional  modules'. It recovers the prescribed multiplicities
but furthermore allows to consider a tensor product of `finite
dimensional  modules' realised by a multiplication of their
characters. The result is a new commutative ring which extends
the  ring of $\bgo$ characters. It will be used in sections 6,7
as a basis for the construction of quantised  `$q$'-characters
which realise the fusion rules of the admissible representations 
at rational level.

We can describe the analogs of the multiplicities in \cch\ via
the map $\iii$, so what remains to be done is to generalise the
formal exponentials entering \cch.  The naive extension of the
map $\iii$ to the characters, thus leading to ordinary formal
exponentials with arguments involving the weights $\iii(y)$ (or
the horizontal projections $\bL=h(y)
\in \bh^*$) is possible but is not consistent (except in the case
$\bgo=sl(2)$) since in general a sum of such weights  goes beyond
the image of $\iii$. (This is in agreement with the fact that
there is no  nontrivial subset of the $n$-ality zero
representations of $\bgo$ closed under tensor products).  So our
main idea is, instead of  exponentials (elements of the group
algebra of the group of translations $t_P$) assigned to  weights,
to consider the elements $ \bw\,\ts-,{\zl}=y\in \tW\,,$ i.e., the
generating elements of the group algebra of the group $\tW$ (the
extension of $t_P$ over $\bW$). (We shall use the same notation
for both interpretations denoting as before sometimes the affine
translations by formal exponentials.)\foot{ 
{}From now on $\hw\, \ts,\zl$ will stand for the
multiplication of two elements in the group ring $\IZ[\tW]$ ;
accordingly the standard notation for the  action of the
(horizontal) Weyl group on the formal exponentials,
$w(\ts-,{\zl})=\ts-,{\hw(\zl)}\,, $ will be replaced by
$\ts-,{\hw(\zl)}=\hw\ts-,\zl\hw^{-1}$, cf. \Ie.  }
Unlike the formal exponentials these generating elements 
do not commute any more 
(rather satisfy the multiplication rules of $\tW$ in \Ie,
but nevertheless  the `finite dimensional module' characters we
obtain, do commute and multiply according to the rules
conjectured in our previous paper.
\medskip

Now we introduce the supports, $\VV_{u}\,$ and  $\,\GG_{u}$, and 
characters, $\ch(\VV_{u})\,$ and $\cc_u$, of, respectively,
`Verma' and `finite dimensional modules'.  The supports are
certain infinite or finite subsets of the extended affine Weyl
group $\tW$ while the characters are certain (formal) series or
finite sums with integer coefficients of elements of $\tW$.

Let $\zL=y\cdot k\zL_0\,,\  y\in a W \,, a\in A\,.$ We define
\eqn\vsupp{
\VV_{y}:=\{z\in a W |\, \iii(z)\in \iii(y) - Q^+   \}
=\{xy\,|\, \ x\in W \,, \ \iii(x)\in  - \ol{y}(Q^+)   \}\,.
}
Alternatively, using that any $\nu\in Q^+$ can be represented
uniquely as $\nu=n\beta+\zl\,,$ with some $\zb\in Q^+$ and
$\zl\in Q^+/ n Q^+=\{\sum_{i=1}^{n-1}\, k_i \za_i\,|\, 0\le
k_i\le n-1\}$, we can bring $\VV_y$ into the form
\eqn\vsuppa{
\VV_{y}= \HT^{\ol{y}}\, y \, t_{Q^+}=\{ t_{_{\ol{u} \ol{y}(\zb)}}
\, u y \,|\,u\in \HT^{\ol{y}}\,,\  \zb \in Q^+ \}\,.
}
Here the finite subsets $\HT^{\ol{y}}\subset W$, which project
horizontally to $\bW\,,$ are subject to the condition
\eqn\ty{
  \HT^{\ol{y}}=\{u\in W\,|\, -\by^{-1}(\iii(u))\in Q^+/ n Q^+\}\,.
}
E.g., if an element $u\in \HT^{\ol{y}}$ projects to a reflection
$\bw_{\za}\in \bW\,, \za\in \bR_+\,$ then it equals
$u=\bw_{\za}\,,$ or $u=\bw_{\za}\, t_{-\za}=w_{\zd-\za}$ if
$\ol{y}^{-1}(\za)\in \bR_+\,,$ or $\ol{y}^{-1}(-\za)\in \bR_+\,,$
respectively.

Thus the support $\VV_y$ in \vsupp\ naturally generalises the
support of a $\bgo$ Verma module,  being defined as a collection
of $|\bW|$ ``positive'' (for a choice of a set of simple roots)
root lattice cones $t_{\ol{u}\,\ol{y}(Q^+)}\,$ applied to $ u
\,y\in \HT^{\ol{y}}\,y\,.$ The two descriptions of $\VV_y$ 
reflected in \vsuppa\  and in \vsupp\ -- as a collection of
supports of ordinary $\bgo$ modules of highest weights $h(u\,y
)=\overline{u y \cdot k\zL_0}$, or as  the support of $\bgo$
module of highest weight $\iii(y)$ with  some `excluded' points,
will be both useful in what follows.  Next we define the
multiplicity of an weight $z$ through its $\iii$ image,
\eqn\vmul{
\KK^{y}_{z} :
=K_{\iii(y)-\iii(z)} \, \quad {\rm or}\ \  
\KK^{y}_{xy}  = \KK_{-\ol{y}^{-1}(\iii(x))} \,,
}
where $K_{\zb}$ is the Kostant partition function in \cch.
Then the characters generalising \cch\ are defined according to

\eqn\vch{
\ch(\VV_{y}):
 =\sum_{\matrix{z\in \tW\,, zy^{-1}\in W\,\cr \iii(z)\in
\iii(y)-Q^+}} z\, \KK^{y}_{z}   \,
 = \sum_{\matrix{x\in W\,,\cr \iii(x)\in -\ol{y}(Q^+)}} xy\,
\KK_{-\ol{y}^{-1}(\iii(x))}\,.
}
Equivalently, using that
$-\by^{-1}(\iii(x))=-\by^{-1}(\iii(u))+n\zb$ for $xy=u y
t_{\zb}\,,$ and denoting
\eqn\pk{
P_{r}\, \ts-,{\br} \,\dd: =\sum_{\zb\in Q^+}\, K_{n \zb+r}\,{\tt
e}^{ {\kappa\,\zb} } \,,
}
we rewrite \vch\ as
\eqn\vcha{
\ch(\VV_{y})
= \sum_{u\in\HT^{\ol{y}}} u y\, 
P_{-\by^{-1}(\iii(u))} \, \ts-,{\br}\,\dd\,.
}
Multiplying both sides of \pk\ by  ${\tt e}^{ {\kappa\over n}
\,r}\,,$ summing over $r\in Q^+/n Q^+$ and making a change of
variables $n\zb+r=\zb'\,,$  we recover in the r.h.s.  the factor
${\tt e}^{- {\kappa\over n} \br}\,d_{\kappa/n}$, cf. \cch.  Thus
$P_r$ are
polynomials of translations in 
$Q^+$  
determined through 
the quotient 
of the standard denominators, $d_{\kappa/n}$ and
$\dd$.
\eqn\pI{
\sum_{r\in Q^+/n Q^+}\, P_{r}\,{\tt e}^{ {\kappa\over n} \,r}=
{{\tt e}^{- {\kappa\over n} \br}\,d_{\kappa/n}\over  \ts-,{\br}
\,\dd}  = \prod_{\za>0}\, \sum_{k_{\za}=0}^{n-1}\, {\tt e}^{
{\kappa\over n}\,  k_{\za}\,\za }=\sum_{\mu\in  Q^+
 }\, {}^{(n)}K_{\mu}\, {\tt e}^{{\kappa \,\mu\over n}}\,,
}
\eqn\pII{
P_r=
\sum_{\nu\in Q^+}\, {}^{(n)}K_{r+n\nu}\, {\tt e}^{
{\kappa \,\nu}}\,.
}
 The partition function  ${}^{(n)}K_{\mu}\,$ defined through the
last equality in \pI\ is apparently nonzero for a finite subset
of $Q^+$, i.e., the summation in \pII\ is finite.  The relations
\pk, \pII\ imply
\eqn\pIII{
K_{n \zb+r}\,=
 \sum_{\zn\in Q^+}\, {}^{(n)}K_{n\nu+r}\, \,
  K_{\zb-\zn}\,  \,.
 }

We have  the symmetry properties 
\eqn\spr{\eqalign{
 \ch(\VV_{_{a\, y}})&= a\,
\ch(\VV_{_{y}})\,, \cr
a^{-1}\, \HT^{\ol{a}\ol{y}}\, a &=\, \HT^{\ol{y}}
\,,\quad
P_{-({\ol{a}\ol{y}})^{-1}(\iii(aua^{-1}))}=P_{-\by^{-1}(\iii(u))}
\,,
}}
using that $\iii(a)=0\,,$  $\iii(a^{-1} x
a)=\bar{a}^{-1}(\iii(x))$ for $a\in A\,.$

Now consider  $\zL=y\cdot k\zL_0\in \CP_+\,,$  $y\in \tC\,,$
hence $\iii(y)\in P_+$.  According to our definition $\VV_{y}$ is
reducible with submodules  $\VV_{y\bw}\,,$ $\bw\in \bW\,,$ since
$V_{\iii(y)}$ is reducible with submodules
$V_{\bw^{-1}\cdot\iii(y)}\,.$ In parallel with \sfch\ we define
the characters of `finite dimensional modules' by a  `resolution'
formula with respect to the KW Weyl group
\eqn\res{
\cc_{y}\,:  =  \sum_{\hw'\in\bKW{\zL}} \det(\bar{w'})\
\ch(\VV_{w'\,y}) = \sum_{\bw\in\bW} \det(\bw)\ \ch(\VV_{y\,\bw}) 
            =  \sum_{z\in\tW\,,\, zy^{-1}\in W} \ml_{z}^{y} \ z\,.
}
We extend this definition to the whole $\tW$, i.e.,
\eqn\resa{
\cc_{y\bw}=\det(\bw)\,\cc_{y}\,, \quad y\in \tC\,, \ \bw\in
\bW\,. 
}
Using \vmul\ and the intertwining property \Ic\
of the map $\iii$ \res\ gives for the multiplicities $\ml_{z}^{y}$ 
(cf. \sfch\ )
\eqn\mt{
\ml_{z}^{y}=\sum_{\bw\in\bW} \det(\bw)\, \ K_{z}^{y \bw}=
\sum_{\bw\in\bW} \det(\bw)\, \ K_{\iii( y \bw)-\iii(z)}
=\sum_{\bw\in\bW} \det(\bw)\,
  \ K_{\bw\cdot\iii(y)-\iii(z)}
=\bm^{\iii(y)}_{\iii(z)}\,.
}

Having an explicit description for the multiplicities we can
introduce the supports $\GG_{y}$ of `finite dimensional modules'
as $\GG_{y}=\{z\in \tW \,|\, \ml_{z}^{y}\ne 0\}$.  {}From \mt\
and from the definition of the map $\iii$ it follows that these
`generalised  weight diagrams' have the structure of  $n$-ality
zero $sl(n)$ weight diagrams $\zG_{\iii(y)}$ with the  points
$\mu\not\in\ $Im$(\iii)$ excluded.  An $sl(3)$  example is
illustrated on fig. 4, see \FGP\ for more examples.   Let us
point out some symmetry properties of the
generalised  weight diagrams. We have
from \spr
\eqn\scur{ 
\cc_{_{a\,y}}=a\,\cc_{_{y}}\,,   \quad a\in A\,,
}
 which implies  that $m^{a y}_{a z}=m^y_z$.
The invariance of the $\bgo$ multiplicities
$ \bm_{\bar{a}(\iii(z))}^{\iii(y)}= \bm_{\iii(z)}^{\iii(y)}=
 \bm_{\iii(a\,z\,a^{-1})}^{\iii(y)}\,$) implies
\eqn\mti{
m_{a^{-1}\,z\,a}^y=m_z^y\,,\quad
 a\in A\,
}
and hence
$a\,\cc_{_{y}}\,a^{-1} =\cc_{_{y}}\,.$
  For $z= a\,
t_{-\mu}\,,$ $a\in A\,,$ the symmetry property \mti\ extends to
\eqn\mtia{
m_{a\,t_{-\bw(\mu)}\,}^y=m_{a\,t_{-\mu}}^y\,,\quad \bw\in \bW\,,
}
using once again \tlg\ to obtain 
$\iii(a\bw\,a^{-1}\,z\,\bw^{-1})
=\bw(\iii(z))$ as well as the invariance
of the $\bgo$ multiplicities
$ \bm_{\bw(\iii(z))}^{\iii(y)}= \bm_{\iii(z)}^{\iii(y)}$ for any
$\bw\in \bW$.

 \medskip
Now we examine the case $\bgo=sl(3)$ in details.  The  6 element
sets $\HT^{\ol{y}}\,$ are $\HT^{\bar{a}}=a\,\bW\,a^{-1}\,,$ while
$\HT^{\bar{a} w_{\theta}}=a\,\HT^{w_{\theta}}\,a^{-1}\,,$ and
$\HT^{w_{\theta}}= \{\un\,,
w_{010}\,,w_{020}\,,w_{10}\,,w_{20}\,, w_0\}\,,$ see figs. 5,6,
where the two basic supports are  `visualized' as certain
subgraphs of the Cayley graph of $W\,.$ The set of polynomials
$P_{-\by^{-1}(\iii(u))}$ in \vcha\ associated with $\HT^{\ol{y}}
$ and the corresponding characters for $\ol{y}=\un\,, w_{\theta}$
read
\eqn\cha{\eqalign{
&   \ch(\VV_{_{t_{-\zl}}}) = 
 \sum_{\bw\in \bW} \bw\, \, 
 P_{-\iii(\bw)} \, \ts-,{(\zl+\br)}\
\dd \cr
&=
    \bigl[  (1+2\ts,{\theta}  ) 
          + \hw_{12}(2+\ts,{\za_1}) 
          + \hw_{21}(2+\ts,{\za_2}) 
\cr &
                  + \hw_1(1+\ts,{\theta}\,+\ts,{\za_2}) 
          + \hw_2(1+\ts,{\theta}\,+\ts,{\za_1})
          + 3\hw_{\theta} 
   \bigr] \ts-,{(\zl+\br)} \ \dd \,,
}}
\eqna\chaa
$$\eqalignno{
&   \ch(\VV_{_{w_{\theta}\,t_{-\zl}}}) = 
     \sum_{u\in \HT^{w_{\theta}}} u\,w_{\theta} \,
         P_{-\theta^{-1}(\iii(u))}\,\, \ts-,{(\zl+\br)}\ 
\dd\cr
& =
        \bigl[    \hw_{0}\,\hw_{\theta}\,    (2 + \ts ,{\theta}) 
            + \hw_{010}\,\hw_{\theta}\, (1 +2\ts ,{\za_2}) 
            + \hw_{020}\,\hw_{\theta}\, (1 +2\ts ,{\za_1})
&\chaa {}\cr &
            + \hw_{20}\,\hw_{\theta}\,  (2 + \ts ,{\za_1}) 
            + \hw_{10}\,\hw_{\theta}\,  (2 + \ts ,{\za_2}) 
                    + \hw_{\theta}\,  (1 +2\ts ,{\theta}) 
                 \bigr]\, 
   \ts-,{(\zl+ \br)} \ \dd \,.
}$$

\medskip

Now we rewrite the `finite dimensional module' characters
\res\ expressing them in terms of the ordinary $sl(3)$
characters, the elements of $A$ and the `class' element
\eqn\df{
 \FF \equiv \sum_{\aa\in\AA}\, a\,\hw_0\, a^{-1}
       = \hw_0 +  \hw_1 + \hw_2 \ .
}
\PROP {\bf 3.1\ }  {\it For any $ x = \ol{x}
\,t_{-\zl}=t_{-\zn}\,\ol{x}  \in \tC\,, \ \zn=\ol{x}(\zl)$ we
have } 
\eqn\fch{\eqalign{
\cc_{x} 
  & =\det(\ol{x})\, \Big( \bc_{\zn} 
     + \zg\,\bc_{\zn-2\fw_{1}} 
     + \zg^{-1}\,\bc_{\zn-2\fw_{2}} 
   + (\FF+2) \big(\bc_{\zn-\br} 
     + \zg\,\bc_{\zn-\fw_{2}} 
     + \zg^{-1}\,\bc_{\zn-\fw_{1}}  \big)\Big)\cr
&= \bc_{\zl+\ol{x}^{-1}\cdot 0} 
     + \zg\,\bc_{\zl+\ol{x}^{-1}\cdot(-2\fw_{1})} 
     + \zg^{-1}\,\bc_{\zl+\ol{x}^{-1}\cdot(-2\fw_{2})}\cr 
   &+ (\FF+2) \big(\bc_{\zl+\ol{x}^{-1}\cdot(-\br)} 
     + \zg\,\bc_{\zl+\ol{x}^{-1}\cdot(-\fw_{2})} 
     + \zg^{-1}\,\bc_{\zl+\ol{x}^{-1}\cdot(-\fw_{1})}  \big)\,.
}}
\medskip
\PROOF 
First the  `Verma module'  characters \cha, \chaa{} can be
rewritten for any $x= \ol{x} \,t_{-\zm} \in \tW$  as (recall that
$\zg=w_{12}\,\ts-,{\fw_2}$)
\eqn\chp{\eqalign{
 \ch(\VV_{x})& = 
\bigl[ ( 1  + \gg     \ts ,{\ol{x}^{-1}(2\fw_1)} 
                        + \gg^{-1}\ts ,{\ol{x}^{-1}(2\fw_2)} )\cr
&  + (\FF+2) (  \ts,{\ol{x}^{-1}(\br)} 
   + \gg     \ts,{\ol{x}^{-1}(\fw_2)} 
   + \gg^{-1}\ts,{\ol{x}^{-1}(\fw_1)}  ) 
   \bigr] \dd \ts-,{(\zm+\ol{x}^{-1}(\br))}\,.
}}
Next apply the `resolution' formula \res\ in the form
$$
\cc_x=
\sum_{w'\in \bKW{x} }\, \det(\bar{w'})\, \ch(\VV_{w'\,x})=
\sum_{\bw\in \bW}\,\det(\bw)\, \ch(\VV_{\ol{x}
\bw\,t_{-\bw^{-1}(\zl)}}) 
$$
$$
=\det(\ol{x})\, \sum_{\bw\in \bW}\, 
\det(\bw)\,\ch(\VV_{\bw\,t_{-\bw^{-1}(\ol{x}(\zl))}})\,.
$$
To obtain the second equality in \fch, i.e., to
express 
the characters 
in terms of $sl(3)$ characters $\bc_{\zl}$ with
$\zl\in P_+\,,$ 
we have used \winv,
so that each term $\det(\ol{x})\, \bc_{\zn-\zn_a}$ turns into
$\bc_{\zl+\ol{x}^{-1}\cdot(-\zn_a)}\,.$ \endPROOF
Taking into account the same property \winv\ (which in particular
implies that
$\bc_{\zl} =0$ for weights on the shifted reflecting hyperplanes,
i.e., $\bw_{\za}\cdot\zl =\zl$ for some $\za\in \bR_+$)
some terms in the general expression \fch\ may vanish
 or compensate each other
for some `boundary' weights.

The character formula \fch\ can be also rewritten as a
decomposition over $\bW\,,$ or $\UU$
$
 \cc_{y} = \sum_{\bw\in\bW} \bw \ B^{(\bw)}_{y} 
           = \sum_{u\in\UU} u \ C^{(u)}_{y} 
$
with some coefficients 
$ B^{(\bw)}_{y}$ 
($ C^{(u)}_{y}$)
given by linear
combinations of $sl(3)$ characters and formal exponentials.
{}From this decomposition the multiplicity $m_z^y$ in \mt\
can be rewritten as a sum of $sl(3)$ multiplicities.
E.g., rewrite \fch\
 for the character $\cc_x\,,\ $  
\eqna\dech
$$\eqalignno{
&\cc_x=\sum_{\mu\in P}\, P_{\mu}^{\zl}\,\ts-,{\mu}:
=
\sum_{a\in A\,, \mu\in P}\,
\Big((\bm_{\mu}^{\zl+\ol{x}^{-1}\cdot(-\zn_a)}+2 
\bm_{\mu}^{\zl+\ol{x}^{-1}\cdot(-\zn_a^{'})})\,a + 
\bm_{\mu}^{\zl+\ol{x}^{-1}\cdot(-\zn_a^{'})}\, F\,a \Big)
\,\ts-,{\mu}
\cr
&=
\sum_{a\in A\,, \mu\in P}\,
\Big(\big(\bm_{\mu}^{\zl+\ol{x}^{-1}\cdot(-\zn_a)}+2 
\bm_{\mu}^{\zl+\ol{x}^{-1}\cdot(-\zn_a^{'})}\big)\,a +
\big(\sum_{b\in A}\,\bm_{\mu+\mu_{b,a}}^{\zl+
\ol{x}^{-1}\cdot(-\zn_b^{'})}\big)\, 
a\,w_0 \Big) \,\ts-,{\mu}\cr
&= \sum_{\matrix{z=a\,t_{-\mu}\cr
a\in A\,, \mu\in P}} m_z^x\,a\,\ts-,{\mu}
+ \sum_{\matrix{z=a\, w_0\, t_{-\mu}
\cr
a\in A\,, \mu\in P}} m_z^x\ \
a\,w_0  \,\ts-,{\mu}&\dech{}
}$$
where the weights $\zn_a\,, \ \zn_a^{'} \in P$ can be read from
\fch, and the polynomials $P_{\mu}^{\zl}=P_{\bw(\mu)}^{\zl}$ are
invariant under the nonshifted action of $\bW$. The weights
$\mu_{b,a}\in P_+$ result from the resummation of the `odd' term
in \fch, setting $w_j\,a =bw_0 t_{-\mu_{a,b}}\,, b\in A\,$ 
 (since $\tW={\cal U} \ t_{P}$), which
allows comparing the second and the third lines in \dech{} to get
an expression for the multiplicities $m_z^x$ in terms of the
$sl(3)$ ones.
\medskip

Let us give some examples
\eqn\fc{\eqalign{
\cc_{_{a}}&=a\,, \ \  a\in A\,,\cr
\cc_{_{\hw_0}}    &=2 +  \hw_0 + \hw_1 + \hw_2   =2 + F\,, 
\cr
 \cc_{_{f}}&\equiv
 \cc_{_{t_{-\bL_1}}} ={\bc}_{_{\bL_1}}+(1 + F )\, \zg^{-1}
 =\cc_{_{\zg^{-1}}} \, \cc_{_{\hw_{20}}} 
 \cr 
&= \zg^{-1}\, (\hw_{20} + \hw_{12} + \hw_{01} + 
        \hw_0 + \hw_1 + \hw_2  + 1)\,, \cr
\cc_{_{f^*}}&\equiv
\cc_{_{t_{-\bL_2}}} ={\bc}_{_{\bL_2}}+(1 + F )\,\zg
= \cc_{_{\zg}} \,\cc_{_{\hw_{10}}} 
\cr
&= \zg\, (\hw_{10} + \hw_{21} + \hw_{02} + 
        \hw_0 + \hw_1 + \hw_2  + 1) \,.
}}

Finally we define dimensions of the `finite dimensional modules'
by setting  in \res\ every (generating) element $\bw\,\ts-,{\mu}$
in $\IZ[\tW]$ to 1, i.e., for $y=\by \, t_{-\zl}\in \tC$
\eqn\dII{
D_{y} = \sum_{z\in\tW} \ml^{y}_{z}\,.
}
Using the  decomposition formula
\fch\   the dimension can be expressed as a sum  of dimensions
$\overline{D}_{\zl}\,$ of $sl(3)$  finite dimensional
representations -- the final expression can be cast into the form
\eqn\dIII{\eqalign{
D_{y}&=9\,\prod_{\za >0}\, \la \by(\zl)+{\br\over
3}\,,\, \by(\za) \ra+ {{\rm det} (\by)\over 3}\cr
&={1\over 3}\, \,\prod_{\za >0}\, \la \iii(y)+\br
\,,\, \za \ra+ {\det(\by)\over 3}
={2\over 3}\, \overline{D}_{\iii(y)} +  {\det(\by)\over 3}\,.
}}
Since 
$\iii(a)=0\,, \quad a\in A$ we have $D_y=D_{a\,y}\,.$
Here are the first few numbers produced from \dIII: 
$1\,,\,5\,,\,7\,,\,19\,,\,23\,,\,43\,,\,83\,,\,103\,,\dots .$ 
\medskip

\REMARK  In the case $\bgo=sl(2)$ we have 
$\iii(W)\equiv Q=2
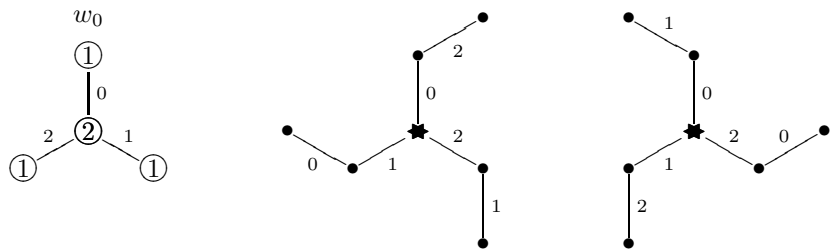
P$ and $\tC\equiv A\,t_{-P_+}$.  The supports of the `Verma
modules' and the weight diagrams ${\cal G}_y$ are isomorphic to
the supports of $\bgo=sl(2)$ Verma modules $V_{\iii(y)}$ and
weight diagrams $\Gamma_{\iii(y)}$, respectively. 
 Alternatively the map 
${1\over 2} \iii$ recovers the supports of
modules of highest weight ${1\over2} \iii(a\,t_{-\zl})=\zl\in
P_+$ of the superalgebra $osp(1|2)$, with
 $a=\un\,, \zg\,,$ labelling the two types of 
modules for a given $\zl\,.$ One has
$\cc_{a\,t_{-\zl}} =a(\bc_{\zl} +\zg\, \bc_{\zl-\za/2}$) and we
can replace $\zg$ simply by a sign $\epsilon=- 1$, thus
recovering the supercharacters of the finite dimensional
representations of $osp(1|2)$.  
The formula \dIII\  is 
replaced by $D_y= \la \iii(y)+\br\,,\, \za \ra=
2\,\la \zl+{\br\over 2}\,,\, \za \ra\,,$ $y=a\,t_{-\zl}\,.$

%
%

\newsec{Character ring and its structure constants}

Denote by $\tCR$, ring of `characters of finite dimensional
modules',  the subring of $\IZ[\tW]$ generated by 
the set of characters $\{\cc_y | y\in\tC\}$ with $\cc_y$
defined in \res.  Thus the multiplication in $\tCR$ is inherited
from the multiplication in the group ring $\IZ[\tW]$.  The main
aim of this section is to obtain a formula for the structure
constants of the ring.

\PROP {\bf 4.1. \ } {\it The ring $\tCR$ 
is a commutative ring.}
\PROOF 
The statement follows from the fact that $A$ is commutative and
commutes with $\FF$, and \ginv\ is equivalent to
$\hw\bc_{\zl}\hw^{-1}=\bc_{\zl}$, i.e., the ordinary characters
commute with elements of $\tW$.
\endPROOF

Next we have

\LEMMA {\bf 4.2 \ } 
\eqn\fmult{\eqalign{
\bc_{_{\fw_1}}\, \cc_{_{y}}& =\sum_{i=1,2,3}\,
\cc_{_{t_{-e_i\,}y}}\,, \quad 
\bc_{_{\fw_2}}\, \cc_{_{y}} =\sum_{i=1,2,3}\,
\cc_{_{t_{e_i\,}y}}\,, \quad 
F\, \cc_y =\sum_{j=0,1,2}\, \cc_{_{{w_j\,y}}}\,.
}}
\PROOF
The first two of these equalities (in which $\sum_i\, e_i=0\,,
e_1=\fw_1\,, e_3=-\fw_2$) follow from the decomposition \fch\ of
$\cc_y$ in $sl(3)$ characters and the analogous multiplication
rules of the latter with the characters of the fundamental
representations $\fw_i\,;$ recall that
$\bc_{_{\fw_1}}=\sum_{i}\,{\tt e}^{-\kappa \, e_i}\,.$ The
derivation of the third is based on the straightforward relation
\eqn\basr{
F^2=3+ \zg\, \bc_{_{\bL_1}}+\zg^{-1}\, \bc_{_{\bL_2}}\,,
}
and the use of \scur\ and the first two of the equalities in
\fmult, taking into account the splitting of \fch\ into `even'
and `odd' part, $\cc_y=\cc_y^{(+)} +F\, \cc_y^{(-)}\,,$
$\cc_y^{(\pm)}$ being linear combinations of $A$ with  $sl(3)$
characters as coefficients (or, elements in the group ring
$\overline{\CR}[A]$ of $A$, over the ring $\overline{\CR}$ of
$sl(3)$ characters). 
\endPROOF

 The quantities in  \fmult\ are the basic
ingredients of the examples \fc\ in the previous section and thus
from the above Lemma and from \scur\ we have

\COR {\bf 4.3}
\eqna\fmt
$$\eqalignno{ 
       \cc_{\zg}\ \cc_y&= \cc_{\zg y}\,,   \cr
  \cc_{\hw_0} \ \cc_{y}&=2 \cc_{y} + \cc_{\hw_0\, y} +
  \cc_{\hw_1\, y}+    \cc_{\hw_2\, y}=\sum_{ x \in \GG_{w_0} }\,
  m_x^{w_0}\,  \cc_{x\, y}= \sum_{ z \in \GG_{w_0}\,y}\,
  m_{z\,y^{-1} }^{w_0}\,  \cc_{z} \,,\cr 
  \cc_{f}\,  \ \cc_{y} &=   \cc_{ f\,y} +\cc_{w_1\, f\,y}+
  \cc_{w_{21}\, f\,y}+\cc_{w_{121}\, f\,y} +\cc_{w_{0121}\, f\,y}
  +\cc_{w_{021}\, f\,y}   +\cc_{w_{2021}\, f\,y}\,, \cr 
                       &=\sum_{ z \in \GG_{f} }\,  \cc_{z\, y}=
  \sum_{ z \in \GG_{f}\,y }\,  \cc_{z}\,,\qquad
  f=t_{-\fw_1}\,,
  &\fmt {}  \cr  
  \cc_{f^*}\,\ \cc_{y} &=   \cc_{ f^*\,y} +\cc_{w_2\, f^*\,y}+
  \cc_{w_{12}\, f^*\,y}+\cc_{w_{212}\, f^*\,y}  +\cc_{w_{0212}\,
  f^*\,y} +\cc_{w_{012}\, f^*\,y}   +\cc_{w_{1012}\,
  f^*\,y}\,,\cr 
                       &=\sum_{ z \in \GG_{f^*} }\,  \cc_{z\, y}=
  \sum_{ z \in \GG_{f^*}\,y }\,\cc_{z}\,,\qquad 
  f^*=t_{-\fw_2}\,.
}$$ 
These Pieri type formul\ae\ hold for generic $y$, in general
there could be cancellations on KW orbits due to \resa.

In \fmt{} the shifted weight diagram $\GG_{f}\,y$, consisting
generically of 7 points, appears, thus we recover the
multiplication rule of the `fundamental' representation $f$
obtained in \FGP.  We recall that it was found  solving the
null-decoupling equations resulting from a pair of singular
vectors of weight $w_{f(\za_i)}\cdot \zL\,,$ $i=1,2\,$ (i.e.,
$w_{\zd+\za_1}\cdot\zL\,,$ $\ w_{\za_2}\cdot \zL\,$)  in the
$\go$ Verma module of highest weight $\zL =t_{-\fw_1}\,\cdot
k\zL_0$. Similarly the general property \scur\ (first line in
\fmt{}) was confirmed in \FGP\ in the case $sl(3)$ analysing the
decoupling conditions corresponding to singular vectors of weight
$w_{\zg(\za_i)}\cdot \zL\,,$ i.e.,  $w_{0}\cdot
\zL\,,$  $w_{2}\cdot \zL\,,$ in the $\go$ Verma module of h.w.
$\zL=\zg\cdot k\zL_0$.

The second of the product rules in \fmt{}  appeared in \FGP\
as a consequence of a conjectured general Weyl-Steinberg type 
formula for the structure constants. We shall now prove that this
formula indeed holds for the multiplication of the characters
constructed in the previous section.  
We proceed in full analogy with
the proof of its standard $sl(3)$ analog. First we establish

\LEMMA {\bf 4.4} {\it For any $x,y\in \CC$ we have}
\eqn\weylI{
     \cc_x\ \cc_y = \sum_{z\in W} \ml^{x}_{z} \ \cc_{zy}
  =\sum_{z\in W\cap \CG_x\, y} \ml^{x}_{zy^{-1}} \ \cc_{z} \,.
}

The formula \weylI\ is the analog of the $sl(3)$ geometric
algorithm for the decomposition of the  product, namely one takes
the generalized `weight diagram' $\GG_x$ and `translates' it by
$y$.  The proof of Lemma 4.4 relies on the decomposition \fch\ of
the characters $\cc_y$ in terms of ordinary characters $\bc$.
(Below we will give an alternative proof.) Recall that the proof
of the  standard analogs of \weylI\ is based on the invariance of
the multiplicities $\bm^\zn_{\bw(\zm)}=\bm^\zn_{\zm}$ \ginv\
under the ordinary (unshifted) action of $\bW$, so that
$\bc_{\zl}=\sum_{\mu\in P_+}\, \bar{m}^{\zl}_{\mu}\,
\sum_{\bw\in \bW}\,\ts-,{\bw(\zm)}\,,$   and the property
\eqn\IVa{
 \sum_{\bw\in\bW}\ \ts-,{\bw(\zm)} \, \bc_{\zl} 
= \sum_{\bw\in\bW}\ \bc_{\bw(\zm)+\zl} \,.
}
Let us now  prove Lemma 4.4.
\PROOF  Using that $\sum_{\bw\in\bW}\,
\ts-,{\bw(\mu)}$ commutes with any $w\in \tW$ \IVa\ implies
\eqn\IVb{
 \sum_{\bw\in\bW}\ \ts-,{\bw(\zm)} \, \cc_{y} 
= \sum_{\bw\in\bW}\ \cc_{_{t_{-\bw(\mu)}\,y
}}\,. 
}
Now we start from the decomposition \dech{}\ for $\cc_x$
and we  insert \IVb\ in $\cc_x\, \cc_y$ using furthermore
 \fmt{}, 
\eqna\IVc
$$\eqalignno{
&\cc_x\ \cc_y &\IVc{}\cr
&=\sum_{\mu\in P}\, P_{\mu}^{\zl}\,\ts-,{\mu} \ \cc_y
=\sum_{\mu\in P_+}\,
P_{\mu}^{\zl}\,\sum_{\bw\in\bW}\,\ts-,{\bw(\mu)} \ \cc_y=
  \sum_{\mu\in P_+}\,
P_{\mu}^{\zl}\,\sum_{\bw\in\bW}\, \cc_{_{t_{-\bw(\mu)}\,y}}
=\sum_{\mu\in P}\,
P_{\mu}^{\zl}\, \cc_{_{t_{-\mu}\,y}}
\cr
&
=
\sum_{a\in A\,, \mu\in P}\, \Big(
(\bm_{\mu}^{\zl+\ol{x}^{-1}\cdot(-\zn_a)}+
2 \bm_{\mu}^{\zl+\ol{x}^{-1}\cdot(-\zn_a^{'})})
\cc_{_{a\,t_{-\mu}\,y}} +
\bm_{\mu}^{\zl+\ol{x}^{-1}\cdot(-\zn_a^{'})}\,
\sum_{j=0,1,2}\, 
\cc_{_{w_j a t_{-\mu}\,y}}\Big)\cr
&
=
\sum_{a\in A\,, \mu\in P}\, \Big(
(\bm_{\mu}^{\zl+\ol{x}^{-1}\cdot(-\zn_a)}+
2 \bm_{\mu}^{\zl+\ol{x}^{-1}\cdot(-\zn_a^{'})})
\cc_{_{a\,t_{-\mu}\,y}} +\big(\sum_{b\in A}\,
\bm_{\zm+\zm_{b,a}}^{\zl+ 
\ol{x}^{-1}\cdot(-\zn^{'}_b)}\big)
\cc_{_{ a w_0 t_{-\mu}\,y}}\Big)\,,
}$$
repeating  the  resummation 
in the second line of \dech{}; it remains to use the 
last line in \dech{}\ to recover \weylI.
 \foot{Strictly speaking  the summation  
over $\bW$ in the first line of \IVc{}, which  
is rather a summation over orbits,
should contain an additional
factor  for weights
$\mu$ with nontrivial stationary subgroup of $\bW$; the same
remark applies to the summation in the last line before \IVa. }
\endPROOF
The lemma extends to $x,y\in \tC$ using \scur, \mti. Next we have
 
\PROP {\bf 4.5} {\it Let $x\in a\,\CC\subset \tC\,, y\in a'\,\CC\,, 
a,a'\in A$. Then }
\eqn\weyl{
  \cc_x\ \cc_y 
  = \sum_{z\in a a'\CC} N^{z}_{x,y} \ \cc_{z} \,,
}
\eqn\strc{
N^{z}_{x,y} 
= \sum_{w'\in\bKW{z}} \,\det(\overline{w'})\ \ml^{x}_{w'\,z\,
y^{-1}}  
= \sum_{\bw\in\bW} \, \det(\bw)\ \ml^{x}_{z\bw y^{-1}} 
= \sum_{\bw\in\bW} \,\det(\bw)\ \bm^{\iii(x)}_{\iii(z\bw
y^{-1})}\,. 
}
Note that the summation in \weyl\ runs effectively over the
shifted weight diagram $\tC\cap \CG_x\,y\,$ (of `shifted highest
weight' $xy\,$) since from the expression of the structure
constants $N^{z}_{x,y}$ it follows that $z\bw\in \CG_x\,y$ for
any $\bw\in \bW$.  To make contact with the notation in \FGP,
where we have used the horizontal projections  of the
weights $\zL_y=y\cdot k\zL_0$, note that the `shifted highest
weight' was denoted in \FGP\ by $\ol{\zL}_x\circ\ol{\zL}_y
=h(x)\circ h(y):=h(x)+ \bar{x}(h(y)$ which coincides,
according to \hor{}, with $ h(x y)$.

\PROOF
Using Lemma 4.4  it remains
to account for  cancellations on KW orbits due to \resa, using
that $\tC$ is a fundamental domain in $\tW$ with respect to the
KW Weyl group.
\endPROOF

Finally we prove another property announced in \FGP.  The
structure constants $N^{z}_{x,y}$ of the character ring $\tCR$
can be expressed by  the structure constants of the $sl(3)$
character ring  $\overline{\CR}$ through the map $\iii$:

\PROP {\bf 4.6}
\eqn\strci{
N^{z}_{x,y} = \bar{N}^{\iii(z)}_{\iii(x)\,\iii(y)} \,.
}
\PROOF
Using \tlg\ we have $\iii(z\bw y^{-1}) = \ol{y}(\iii(z\bw)) +
\iii(y^{-1}) =\ol{y}(\iii(z\bw) - \iii(y))\,. $ Hence $
m^{x}_{z\bw y^{-1}}= \bm^{\iii(x)}_{\iii(z\bw y^{-1})} =
\bm^{\iii(x)}_{\ol{y}(\iii(z\bw) - \iii(y))}
= \bm^{\iii(x)}_{\iii(z\bw) - \iii(y)} =
\bm^{\iii(x)}_{\bw^{-1}\cdot\iii(z) - \iii(y)}$ (using \Ic\ in
the last step), which inserted in \strc\ converts it into the
classical Weyl--Steinberg formula for the $sl(3)$  tensor product
multiplicities $\bar{N}^{\iii(z)}_{\iii(x)\,\iii(y)}$.
\endPROOF

\COR {\bf 4.7} {\it The structure constants 
$N_{x,y}^z$ are nonnegative integers.}

We introduce an involution in $\tW$ induced by the $\IZ_2$
automorphism of the $\bgo \ $ Dynkin diagram $\za_i\rightarrow
\za_i^*:= \za_{n-i}\,,$  $i=1,2,\dots,n-1\,,$ according to
$w_{\za_i}^*=w_{\za_i^*}\,,$ $t_{\zl}^*=t_{\zl^*}\,, \
\la\zl^*\,, \za_i\ra= \la\zl\,, \za_i^*\ra\,,$ and then
$\iii(x^*)=(\iii(x))^*$.  The involution extends to an
automorphism of the ring $\tCR$ with $\cc_y^*=\cc_{y^*}$.
Proposition 4.6 implies the standard properties of the structure
constants
\eqn\invo{
N_{x,y}^{\un}=\zd_{x,y^*}\,,\quad
N_{x,y}^{z}=N_{x,z^*}^{y^*}
=N_{x^*,y^*}^{z^*}
\,,
}
along with 
\eqn\invoi{
N_{x,y}^{z}
=N_{a x,y}^{a z}\,, \quad a\in A\,.
}

\REMARK The `classical' dimensions $D_x=D_{x^*}$  \dII\ provide a
numerical realisation of the product rule \weyl\
$$
 D_x\ D_y   = \sum_{z}\, N^{z}_{x,y} \ D_z
=\sum_{z}\, \bar{N}^{\iii(z)}_{\iii(x)\,\iii(y)} \ D_{z} \,.
$$

The ring $\tCR$
has a set of `fundamental' characters that generate it
as a polynomial ring. One possibility for this fundamental set
is $\{ \cc_{\hw_0}, \cc_{f}, \cc_{f^*} , \cc_{\gg}=\gg \}$.
The group $A$ is like a set of simple currents or in the
language of rings it is a group of units of the ring $\tCR$.
Thus 
$\tCR = \CR[A]$
where  $\CR$ is the
triality zero subring $\CR$ of $\tCR$ having as a
fundamental set $\{\cc_{\hw_{0}}\,, \, \cc_{\hw_{20}}
\,,\, \cc_{\hw_{10}}\}
=\{\cc_{\hw_{0}}\,, \,\zg\, \cc_f\,,\, \zg^{-1}\, \cc_{f^*}\} $.
It is convenient to introduce the notation $\ff_0=\cc_{\hw_0}$,
$\ff_1 = \cc_{\hw_{20}}$, $\ff_2 = \cc_{\hw_{10}}$. The weight
diagrams 
$\CG_{f_j}\,, \ j=0,1,2\,$ are depicted on fig. 7. Since
$\bc_{\bL_1}=\zg^{-1} \, \sum_{a\in A}\, a\,\bar{\zg}\,a^{-1}$ 
and $\bc_{\bL_2}= \zg \,\sum_{a\in A}\, a\,\bar{\zg}^{-1}\,
a^{-1}$ the relations \fmult\ read equivalently

$$
\sum_{\aa\in\AA}\ 
a\,x\,a^{-1}\ \cc_{_{y}}
            = \sum_{\aa\in\AA}\ \cc_{_{
                a\,x\,a^{-1}\ y}} 
$$
for any $y\in\tW$ and $x=\hw_{0},\hw_{10},\hw_{20}$.
Similarly the  Pieri  type formul\ae\ 
\fmt{} can be rewritten  for the triality zero counterparts
$\{f_j\,, \, j=0,1,2\} \,.$  With the  shorthand notation 
$M^0_u=\ml^{\hw_{0}}_u$, 
$M^1_u=\ml^{\hw_{20}}_u$, 
$M^2_u=\ml^{\hw_{10}}_u$, $u\in W$, we have 
\eqn\pieri{
  \ff_{j} \ \cc_{y} 
= \sum_{ u \in \GG_{\ff_{j}} } \, M^j_u\, \cc_{u y } 
= \sum_{ z \in \GG_{\ff_{j}}\,y } \, M^j_{zy^{-1}}\,  \cc_{z} \,.
}

Before we proceed with the next proposition we introduce on the
chambers  $\CC$ and $\tC$ a filtration and related gradation
using the reduced length of words $\fCC{k}=\{ x \in\CC \,|\,
\ell(x) \le k \}$ and $\gCC{k}=\{ x \in\CC \,|\, \ell(x)  =  k
\}$.

\LEMMA {\bf 4.8} {\it Let $x\in\gCC{k}$ then the sets
$\GG_{\ff_{i}}\,x$, for $i=1,2$, have a single element in
$\gCC{k+2}$ while the rest are in $\fCC{k+1}$.}
\PROOF The statement is proved by a direct check, cf. figs.1,7, 
using \pieri.
\endPROOF

\PROP {\bf 4.9}
{\it The ring $\CR$ is generated as a polynomial ring by
$\cc_{\hw_0}$, $\cc_{\hw_{10}}$, $\cc_{\hw_{20}}\,,$ subject to
one algebraic relation}
\eqn\rel{
\cc_{\hw_0}\ \cc_{\hw_0}=2\,\cc_{w_0}+
\un+ \cc_{\hw_{10}}+\cc_{\hw_{20}}\,.
}
\PROOF 
We want to show that any $\cc_{y}$, $y\in\CC$, can be represented
as 
\eqn\ring{
 \cc_y = \sum_{\ze=0,1} \sum_{n_1,n_2\in\IZp}
           c^y_{\ze,n_1,n_2} \, 
                   \ff_0^\ze \ff_1^{n_1} \ff_2^{n_2}   \,, \quad
                   c^y_{\ze,n_1,n_2} \in \IZ\,.
}
 First note that
$\fCC{1}=\{\un\,, \, w_0\}\,,$ 
 and $\gCC{2}=\{w_{10}\,,\, w_{20}\}\,.$
Using
Lemma 4.8 and \pieri, \rel\ we see that these polynomials are
determined inductively going up the gradation.
\endPROOF
Due to the relation \rel\
the ring can be generated also by only two fundamental
characters, i.e., either by $\ff_0,\ff_1$ or by $\ff_0,\ff_2$.
\medskip

Using \pieri\ and \ring\
one can give an alternative proof of Lemma 4.4:
\PROOF 
Let us do it for $\CR$. 
By Proposition 4.9 we have
$$\eqalign{
  \cc_x &= \sum_{\ze=0,1} \sum_{n_1,n_2\in\IZp}
           c^x_{\ze,n_1,n_2} \, 
                   \ff_0^\ze \ff_1^{n_1} \ff_2^{n_2}   
\cr  
  &= \sum_{\ze=0,1} \sum_{n_1,n_2\in\IZp}
           c^x_{\ze,n_1,n_2} \,
         \sum_{\{u,\{u_{i,j}\}\}}          
           (M^0_u)^\ze \Big(\prod_{i=1,2}\prod_{j=1}^{n_i}
           M^i_{u_{i,j}}\Big)\,
                   u^\ze  u_{1,1}\dots u_{1,n_1} u_{2,1}\dots 
u_{2,n_2}\,,
}$$
 using \res\ for each $f_j$ with the notation for the  multiplicities
as in \pieri.
Note that all of the above sums are actually finite sums.
Writing $\cc_x = \sum_{z} \ml^x_z \, z$ we get
$$
  \ml^x_z = \sum_{n_1,n_2\in\IZp} \big(
            c^x_{0,n_1,n_2} \,
            \sum_{\{u_{i,j}\}}     
            \prod_{i=1,2}\prod_{j=1}^{n_i}
            M^i_{u_{i,j}}\,
          + 
                    c^x_{1,n_1,n_2} \,
                   \sum_{\{u,\{u_{i,j}\}\}}        
           M^0_u \prod_{i=1,2}\prod_{j=1}^{n_i}
           M^i_{u_{i,j}}\,       \big)\,,
$$         
where the two internal sums are over subsets $\{u_{i,j}\}$
and $\{u,\{u_{i,j}\}\}$ of $W$ such that
$u_{1,1}\dots u_{1,n_1} u_{2,1}\dots u_{2,n_2}=z$
and $u u_{1,1}\dots u_{1,n_1} u_{2,1}\dots u_{2,n_2}=z$
respectively. 
Using repeatedly associativity together with the Pieri 
type formul\ae\  \pieri,
e.g., 
$$
(\sum_{u_{i,j}} M^i_{u_{i,j}} \, u_{i,j}) \cc_{u_{i,j+1}\dots y}
 = 
 \sum_{u_{i,j}} M^i_{u_{i,j}} \,  \cc_{u_{i,j} u_{i,j+1} \dots y}\,,
$$
we are done.
\endPROOF

\medskip
 \REMARK 
The ring  $\tCR$ can be looked as a `quadratic' extension of the
ring of $sl(3)$ characters $\overline{\CR}$ -- a subring of
$\tCR$ according to \fc.  I.e.,  $\tCR$ is   a polynomial ring
with coefficients in  $\tilde{\gV}:=\overline{\CR}[A]$ modulo a 
quadratic relation, 
$\tCR =\tilde{\gV}[F]/\{F^2 -C=0\} \,; \ \  C\in \gV$ is
given by the r.h.s. of \basr.  
Here $\gV \subset \overline{\CR}[A]$ is a triality
zero subring of $\tilde{\gV}$. 
Similarly $\CR = {\gV}[F]/\{F^2 -C=0\} \,.$

%
%

\newsec{ Alcove of admissible weights  at level $k+3=3/p$.}

Beginning with this section 
we will analyse the case of $\go=\widehat{sl}(3)_k$ at rational
levels, namely 
$\zk=k+3=3/p$ with $p\in\IZ_{\ge 2}\bs 3\IZ$. 
This selects a subseries of the general set of admissible weights
of \KW\  which we will describe in more detail.
\foot{This subset is generic since it is expected
that as in the $\widehat{sl}(2)$ case there is an effective
factorisation of the FR multiplicities for the general admissible
representations at $\kappa=p'/p$ into the multiplicities for the 
two subseries --
at $\kappa=n/p$ and the integrable one at $\kappa=p'$,  the former 
represented  by the r.h.s. of \qstrci, 
which extends to arbitrary $p\in \IZ_{\ge 1}$.}

The rationality of $\zk$ has two consequences -- for $y\in \tW$
the map $y\mapsto y\cdot k\zL_0$ for $y\in \tW$ is not injective
and the KW groups are isomorphic to the affine Weyl group $W$.
Let $\Pi^{^{[p]}}=\{\za_0^{[p]}=p\zd-\theta \}\cup\bP$, and
$\Rr_{+,p}= \bR_+\cup\{mp \zd+\za\,, \za\in
\bR\,, m\in \IZ_{>0}\}$. Denote by $W^{[p]}$ the isomorphic to
$W$ subgroup of $W$, generated by the reflections $\{w_{\za}\,, \
\za\in\Pi^{^{[p]}}\}$.  We have $W^{[p]}=t_{pQ} \rtimes \bW$.
The subgroup $A^{[p]}$ of $\tW$ generated by $
\zg_{_{[p]}}:=t_{p\,\bL_1}\,\bar{\zg}=t_{(p-1)\,\bL_1}\,\zg$
keeps invariant the set $\Pi^{^{[p]}}$ and hence $a \, W^{[p]} \,
a^{-1} = W^{[p]}$ for $a\in  A^{[p]}$.  We have $\tW=W\rtimes
A^{[p]}$ and for $\kappa=3/p$
\eqn\pa{
A^{[p]}\cdot
k\zL_0=k\zL_0\,.
}

Let $y\in\tW$ and ${\cal P}=\{\zL=y\cdot k\zL_0 \,|\, y\in\tW\}$.
{}From the Kac-Kazhdan condition and from the analog of \Ia\ with
$\za\in \Pi^{^{[p]}}$ it is clear that if $y(\za)\in\Rr_+$,
$\forall\za\in\Pi^{^{[p]}}$ the reflections $\{w_{y(\za)},
\za\in\Pi^{^{[p]}}\}$ generate a KW group $W^{[\zL]}$ (to be
denoted also $W^{[y]}$) such that its shifted action on ${\cal P}$
 gives the weights of the Verma submodules of $M_\zL$.
As in \Ia\ the shifted action of $W^{[y]}$ on the weights in $\CP$
is intertwined  with the right action of $W^{[p]}$ on $\tW$.
Moreover $M_\zL$ is a maximally reducible Verma module
with infinitely many singular vectors. Hence
we are led to the definition of 
the alcove
of admissible weights  as  $\CP_{+,p} = \tC_p\cdot k\zL_0
=\CC_p\cdot k\zL_0 $ where
\eqn\Idada{\tC_p =\{y\in\tW \,|\, y(\Pi^{[p]})\subset
\Rr_+ \} \qquad {\rm and}\qquad 
\CC_p = \tC_p \cap W 
\,.} 
Denote 
$ 
P_{+,p}^{(w)}=\{\zl\in P_+^{(w)}
\,|\, \la \zl, \theta\ra < p \ {\rm if}\  w(-\theta)<0 \ {\rm
or}\  \la \zl, \theta\ra \le p\ {\rm if}\ w(-\theta)>0\,, w\in
\bW \} \,.
$

\noindent
In particular $P_{+,p}^{(\un)}$ coincides with the integrable
alcove $P_+^{k}$ at level $k=p-1$.  It is easy to see that the
definition of $\tC_p$ is equivalent to
\eqn\alcIIa{\eqalign{ 
\tC_p& 
= \{ y=\ol{y}\,t_{-\zl}\in \tW \,|\, 
 \zl\in P_{+,p}^{(\ol{y})} \} \cr
&=
A\, t_{-P_+^{p-1}}\cup A\, w_0\  t_{-P_+^{p-2}}=
 t_{-P_+^{p-1}}\, A^{[p]} \cup  w_0\
t_{-P_+^{p-2}}\,A^{[p]}\,,\cr 
}}
\eqn\alcII{
\CC_p=\{ y=\ol{y}\,t_{-\zl}\in W \,|\, 
 \zl\in P_{+,p}^{(\ol{y})}\cap Q \}\,,\quad
\tC_p=\cup_{a\in A^{[p]}}\, \CC_p\, a\,.
}
The second equality in \alcIIa, representing the alcove as a
disjoint union of two leaves, parametrised by the two alcoves
$P_+^{p-1}$ and $P_{+}^{p-2}$, takes into account the equivalence
of elements in $\tC_p$ implemented by the right action of the
group  $A^{[p]}$, or, more explicitly,
\eqn\alca{\eqalign{
t_{-\zl}\, \zg_{_{[p]}}^l= \zg^l\, t_{-\zs^{-l}_{[p-1]}(\zl)}&\,,
\quad \zl\in P_+^{p-1}\,,\cr
w_0\,t_{-\zl}\,\zg_{_{[p]}}^l= \zg^{-l}\,w_0\,
t_{-\zs^{-l}_{[p-2]}(\zl)}&  \,, \quad \zl\in P_+^{p-2}\,.
}}
Here $\zs_{[k]}(\zl):=\overline{\zg(\zl+k\zL_0)}=
w_{12}(\zl)+k\,\bL_1$ denotes the automorphism of the alcove
$P_+^{k}$ at integer level $k$ induced by the action of $A$.
Alternatively,  due to \pa,  the admissible alcove is parametrised
 by the
elements of the fundamental domain $\CC_p$ of $W$ (i.e., triality
zero points on any orbit of $A^{[p]}$ in $\tC_p$), as indicated
in \alcII.

In analogy with  Lemma 2.2  one can show that $P=\left(
\cup_{w\in\bW} \, w(P^{(w)}_{+,p}) \right) + pQ$
\quad is a partition and hence one has that  
$\tC_p$, respectively $\CC_p$, is a fundamental domain in $\tW$, 
respectively $W$, for the
right action of $W^{[p]}$. 
Again the
map $\iii$ intertwines the right action of $W^{[p]}$ on $W$ with
the action of the affine Weyl group at level $3p-3$; it is
sufficient to check, taking into account \Ic, that
\eqn\afiii{
\iii(y
w_{p\zd-\theta})=\overline{w_0\cdot(\iii(y)+(3p-3)\zL_0)}\,. 
}
Accordingly $\tC_p$ is represented by a formula analogous to \If,
with $P_+$ replaced by  $P_+^{3p-3}$.

As an example we depict the alcove of admissible weights
$\CP_{+,p}$, for  $p=5$, on fig. 8.  It is parametrised by $\{
wt_{-\zl}\in\tC_p\,|\,w\in\{\unit,w_{\theta}\}\}$ with a circle
or box in case of $w=\unit$ or $w=w_\theta=w_0\, t_{\theta}$
respectively, the numbers inside being the labels of $\zl$.
Equivalently, 
keeping only the triality zero labels $\zl$, the same figure
depicts the alternative representation of the admissible alcove
through the elements of the fundamental domain $\CC_p\,.$ Unlike
\FGP\ the latter choice will be mostly used in what follows.  
Sometimes it will be also useful to work 
with  the full domain $\tC_p\,$  imposing the constraints
implemented by the right action of $A^{[p]}$.  

Define a triality
preserving  order 3 automorphism of $\tW$ (and hence of $W$)
\eqn\afsc{\eqalign{
\sigma_p(x):&=\zg\,x\,\zg_{[p]}^{-p}\,,\quad x\in\tW\,, \cr
 \iii(\zs_p(x))& = \zs_{[3p-3]}^{p}(\iii(x))\,.
}}
Geometrically $\zs_p$ fixes the `middle' point of the alcove
$\CC_p$ ($t_{-{p-1\over 3}\br}\,,$ or $w_0\,t_{-{p-2\over
3}\br}=w_{\theta}\, t_{-{p+1\over 3}\br}\,,$ cf. \alca{}), and
``rotates'' it sending the ``corners'' into one another, i.e., it
behaves like the usual `simple current' automorphism of an
integrable alcove. 

%
%

\newsec{Quantised `${\bf q}$'-characters}

Recall first the integrable case where the
`classical' $\bgo$  characters $\bc_{\zl}\,, $ $\zl\in P_+\,$ are
converted into $\IC$-valued `$q$'-characters, labelled by the
set   $\{\zl\in P_+^{k}\}$ of integrable highest weights at
(positive) integer level $k$.  Essentially one turns the formal
exponentials ${\tt e}^{\zl}\,, \zl\in P\,$ into `true'
exponentials,
\eqn\Vc{
{\tt e}^{\zl}\rightarrow {\tt e}^{\zl}(\mu):=
{e}^{{-2\pi i\over k+n}\la\zl\,, \mu+\br\ra}
\,, \quad \zm\in P \,.
}
 This `quantises' the  `classical' $\bgo$ characters into
`periodic' 
characters, $\bc_{\zl}^{(h)}(\mu)=\bc_{\zl+
h\, Q}^{(h)}(\mu)$, $h=k+n\,,$
 i.e., (skew-)invariant  under the full
affine Weyl group at level $k$, 
\eqn\VIa{
 \det(\bw)\, \,
\bc_{_{\overline{w\cdot(\zl +k\zL_0)}}}^{(h)}(\mu)= 
\bc_{ \zl}^{(h)}(\mu)=\bc_{ \zl}^{(h)}(\overline{w\cdot(\zm
+k\zL_0)})\,, \quad w\in W\,,
}
so that we can restrict the `dual'  set
(the set of $\mu$'s) to the integrable 
alcove $ P_+^{k}$
itself. 
       \foot{Alternatively the `$q$'-characters are
       obtained  restricting the standard group characters 
       to the discrete subset of elements
       $\{$diag$( {e}^{{-2\pi i\over h}\la e_i\,, \mu+\br\ra}\,, 
       i=1,2,\dots, n)\,, \mu\in P_+^{h-n}\}\,$
       in the Cartan subgroup of $SU(n)$. Here $\sum_{i=1}^n\,e_i=0
	   \,, e_i=\bL_i-\bL_{i-1}\,,\, \bL_0=0=\bL_n\,.$} 
The `q'-characters  are given explicitly by a ratio
$\bc_{\zl}^{(h)}(\mu)=S_{\zl\,\zm}^{(h)}/S_{0\,\zm}^{(h)}$ of
matrix elements of the integrable modular matrix
$S_{\zl\,\mu}^{(h)}$, a unitary, symmetric matrix. It is recovered
up to an overall constant by the second equality in \sfch, with
$\kappa=-1$ and exponentials transformed as in \Vc, i.e., \sfch\
turns into  the Kac--Peterson formula \K.  Thus the complex  numbers
$\{\bc_{\zl}^{(h)}(\mu)\,, \ \mu\in P_+^{h-n}\}$ can be
interpreted as eigenvalues of the matrix $N_{\zl}$ of fusion rule
coefficients $(N_{\zl})_{\za}^{\zb}=N_{\zl\,\za}^{\zb}$ of the
integrable WZW conformal models. This relates the Verlinde
formula for $N_{\zl\,\za}^{\zb}$ to the classical Weyl-Steinberg
formula \K, \MW, \FGPa. In what follows we shall also need
$\bc_{\zl}^{(h)}(\mu)\,$ for $\mu$ belonging to some of the
shifted hyperplanes $H_{\za}^{(l h)}:=\{\mu\in \bar{h}^{*} |\,
\la\mu+\br\,,\za\ra=l\,h \}$, $\za\in \bR_+\,$, $ l\in \IZ$.
While the Kac-Peterson formula has no sense on these hyperplanes,
since both the numerator and the denominator vanish, the
characters $\bc_{\zl}^{(h)}(\mu)$ are well defined through the
analog of the last equality in \sfch, or any of the standard
determinant formul\ae\ for the classical $sl(3)$ characters.

Following the analogy with the integrable  case the idea is to
replace the affine Weyl group with the affine KW group at level
$\kappa-3=3/p-3\,,$  i.e., to extend the invariance \resa\ of the
`classical' characters with respect to the right action of the
horizontal Weyl group $\bW$ to invariance  with respect to the
right action of the affine group $W^{[p]}$.  This will lead to
\Iqch\ with the  structure constants given by the conjectured in
\FGP\ formula \ws, which now derives from the `classical' formula
\strc. Finally  inverting \Iqch\ we will recover in section 7 the
Pasquier--Verlinde  type formula \pasv.

Apparently there are two  problems to be solved. We have to find
an analog of the discrete set $\{\mu\in P_+^{k}\}$ and
furthermore the elements of the group algebra of $\tW$  have to
be converted into some $\IC$-valued functions on this set.

 Denote by $E_p$ the `double alcove' region 
\eqn\ExI{\eqalign{
E_p &=\{\zm\in P_+\, | 0\le \la \zm,
\za_i\ra \le p-1\,, \, i=1,2\} \cr
&=P_+^{p-3}\cup 
\big(w_{\theta}(P_{+}^{p-3}) +(p-2)\theta\big)\cup_{\za\in \bR_+}
\big(H_{\za}^{(p)}\cap P_+^{2p-2}\big) \subset P_+^{3 p-3}\,.
}}
This set, which can be also looked as $P_+^{p-1}\cup
\{w_{\theta}(P_{++}^{p+1}) +p \, \theta \}\,,$ 
contains $p^2$ weights,  $|E_p|\equiv |\CP_{+,p}|\,,$ and we
shall argue below that it is the analog for $k=3/p-3$ of the
integrable `dual' set $\{\mu\in P_+^{h-3}\}$, see figs. 9a, 9b,
where $E_p$ is depicted for $p=5$ and $p=4$ (the dotted lines
indicate the hyperplanes $H_\alpha^{(l p)}$).  For $p=2\ $ $E_p$
consists of the alcove $P_+^1$ and the weight $(p-1,p-1)=(1,1)$
and thus  represents the $\IZ_3$ factorisation of the integrable
alcove at level $3p-3$, $P_+^{3p-3}$, obtained after identifying
the points $\zs_{[3p-3]}^l(\zl)$ along an orbit of the $\zs$
automorphism of $P_+^{3p-3}$, including the $\zs$ stable point
$(p-1,p-1)$. For $p>2$ this factorisation leads to a subset of
the alcove $P_+^{3p-3}$ which is of cardinality $|E_p|+
|P_+^{p-3}|> |E_p|$.
 \medskip

We look for a  solution of the invariance condition
\eqn\Vaa{
\cc_{y\,w}(\cdot)=\det(\bw)\, \cc_{y}(\cdot)\,, \quad w\in
W^{[p]}\,, 
}
together with 
\eqn\Ve{
\cc_{y\,a }(\cdot)=\cc_{y}(\cdot)\,, \quad a
\in A^{[p]}\,, \ \ y\in \tC_p\,.
}
Accounting for the invariance of the characters
with respect to the horizontal Weyl
group $\bW$  \resa\ the requirement \Vaa\ reduces to the
periodicity condition
\eqn\Vaaa{
\cc_{y\,t_{p\,\nu}}(\cdot)= \cc_{y}(\cdot)\,, \quad \nu\in Q\,. 
}
The formula \fch\ for the characters $\cc_{y}(\cdot)$ involves
the three basic ingredients -- the elements of the group $A$, the
$sl(3)$ characters $\bc_{\zl}\,,$  and the combination $F$ in
\df, so we have to give meaning to some $\IC$ -  valued
counterparts $\zg(\cdot)\,,\ $ $\bc_{\zl}(\cdot)\,,\ $
$F(\cdot)\,.$ The natural realisation for the generator of the
group $A$ -- isomorphic to the cyclic group $\IZ_3$, reads
\eqn\Vd{
\zg\rightarrow \zg(\mu):=
{e}^{{2\pi\,i\, m\over 3} \zt(\mu)}\,,
\quad m=1,2\,, \ {\rm mod} \ 3\,.
}
The periodicity requirement \Vaaa\  suggests to look for a 
realisation of the $sl(3)$ characters in \fch\ 
in terms of the integrable characters $\bc_{\zl}^{(p)}(\mu)$ at
(shifted) level $p$, determined for $\mu\in P_+^{p-3}\,, $
\eqn\inch{
\bc_{\zl} \rightarrow \bc_{\zl}(\mu):=
\ze_{\zl\,, \mu}\, \bc_{\zl}^{(p)}(\mu)\,,\quad \ze_{\zl\,,
\mu}^3=1\,. 
}
 In
\inch\  we have allowed for an arbitrary overall phase constant
$\ze_{\zl\,,\mu}\,,$ invariant with respect to both indices
under the shifted action of the affine Weyl group.  We can choose
\eqn\ovs{
\ze_{\zl\,, \mu}= {e}^{{-2\pi\,i\, l\over 3} \zt(\mu)\,
\zt(\zl)}\,, \quad l=1,2\,, \ {\rm mod} \ 3\,,
}
which  effectively leads to the realisation of the formal
exponentials as
\eqn\inex{
\ts-,{\zl}\rightarrow 
\ts-,{\zl}(\mu):=
 {e}^{{-2\pi\,i\, l\over 3} \zt(\mu)\, \zt(\zl)}\,
 {e}^{{-2\pi\,i\, \over p} \la \zl, \mu+\br \ra}\,.
}
The need for this phase is dictated by the requirement
\Ve,
which combined with
\scur,  \alca{} reads 
for each of the  parts $\cc_y^{(\pm)}$ in   $\cc_y=\cc_y^{(+)} 
+F\, \cc_y^{(-)}\,$ 
 (treating for the time being $F(\mu)$ as a 
formal variable)
\eqn\Vg{\eqalign{
\zg(\mu)\,
\cc_y^{(\pm)}
(\mu)&=\cc_{t_{-\zs_{[p-1]}(\zl)}}^{(\pm)}
(\mu)\ \big(
=\cc_{t_{-\zs_{[p-3]}(\zl- 2 \bL_2)}}^{(\pm)}
(\mu)  \big) 
\qquad
{\rm for} \ y=t_{-\zl}\,, \cr
  \zg(\mu)\,
\cc_y^{(\pm)}
(\mu)&=\cc_{w_0t_{-\zs_{[p-2]}^2(\zl)}}^{(\pm)}
(\mu)\ \big(= \cc_{w_0t_{-\zs_{[p-3]}^2(\zl-\bL_1)}}^{(\pm)}
(\mu) \big)
\quad 
{\rm for} \ y=w_0\,t_{-\zl}\,.
}}

The above conditions  and the corresponding standard property
 of the integrable `$q$'- characters
\eqn\ssym{
\bc_{\zs_{[p-3]}(\zl)}^{(p)}(\mu)=
{e}^{{2\pi\,i\over
3}\zt(\mu)}\,\bc_{\zl}^{(p)}(\mu)\,,
}
fix the integer $l$ to $l=p$ mod $3$ (using that $p^2-1=0$ mod
$3$), and keeps arbitrary the power $m$ in the phase in \Vd.
Without lack of generality we can choose $m=l=p$ since otherwise
the remaining phases can be absorbed using the analogous to
\ssym\  symmetry  with respect to the index $\mu$,
\eqn\ssyma{
\bc_{\zl}^{(p)}(\zs_{[p-3]}(\mu))=
{e}^{{2\pi\,i\over
3}\zt(\zl)}\,\bc_{\zl}^{(p)}(\mu)\,,
}
thus changing the value of $\mu$ to
$\mu'=\zs_{[p-3]}^{m-p}(\mu)
\in P_+^{p-3}$; we can do this since the three terms in each of
$\cc_y^{(\pm)}$
are described by $sl(3)$ characters of weights of different triality
$\tau=0,1,2\,.$ 

\medskip
Now we turn to the operator $F=w_0+w_1+w_2$. We recall that it
commutes with the elements of $A$ as well as with the $sl(3)$
characters. Preserving the relation \basr\ -- which is the basic
relation used to derive the character ring structure constants,
we see that the square of $F(\cdot)$ can be determined by the
(fundamental) integrable characters, i.e.,
\eqn\sqf{
F^2\rightarrow F^2(\mu):=3+\bc_{\bL_1}^{(p)}(\mu)
+\bc_{\bL_2}^{(p)}(\mu) 
}
for any $\mu\in P$. This determines $F(\mu)$ up to a sign,
$F(\mu)= \ze(\mu)\, \sqrt{F^2(\mu)}\,, \quad \ze(\mu)=\pm 1\,.$
The r.h.s of \sqf\ is equivalently reproduced  by  
\eqn\froot{\eqalign{
F^2(\mu)&= |R(\mu)|^2
\,, \cr 
R(\mu)&=\sum_{\bar{a}\in \bar{A}}\,
e^{-{2\pi \, i \over 3 p} \la \bar{a}(\theta),\mu+\br\ra}=
e^{-{2\pi \, i \over 3 p} \la \theta,\mu+\br\ra}+
e^{{2\pi \, i \over 3 p} \la \za_1,\mu+\br\ra}+
e^{{2\pi \, i \over 3 p} \la \za_2 ,\mu+\br\ra}  \,.
}}

One has the relations 
\eqn\rrela{
i\,3p\sqrt{3}\, S_{0\, \mu}^{(3
p)}=1/d_{\kappa/3}(\mu)=R(\mu)-\bar{R}(\mu)\,,  
}
\eqn\rrel{\eqalign{
&i\,p\sqrt{3}\,S_{0\, \mu}^{(p)}=1/\dd(\mu)=
\sum_{\bar{a}\in \bar{A}}\,{\tt e}^{- \bar{a}(\theta) 
\kappa }(\mu)-\sum_{\bar{a}\in \bar{A}}\,{\tt e}^{\bar{a}(\theta) 
\kappa }(\mu) =(R(\mu))^3-(\bar{R}(\mu))^3\cr
=&
i\,3p\sqrt{3}\, \, S_{0\, \mu}^{(3 p)}\,
\Big(R(\mu)+\bar{R}(\mu)-|R(\mu)| \Big)\,
\Big(R(\mu)+\bar{R}(\mu)+|R(\mu)| \Big)\,.
}}
(Here  $\bar{R}(\mu)$ is the complex conjugation of $R(\mu)\,.$)

It remains to determine the sign of $\ze(\mu)$.  
Since the  parts $\cc_y^{(\pm)}(\mu)$ in $\cc_y=\cc_y^{(+)} 
+F\, \cc_y^{(-)}\,,$   as well as
$F^2(\mu)$, coincide  for $\mu$ and its reflected images
according to \VIa, we can assign $\ze(\mu)=1$ for $\mu\in
P_+^{p-3}$ and $\ze(\mu)=-1$ for $\mu$ sitting on the `mirror'
(with respect to the hyperplane $H_{\theta}^{(p)}$) alcove in
$E_p\,.$   On the intersection of $E_p$ with the
reflection hyperplanes $H_{\za}^{(p)}$ we choose $\ze(\mu)=1$ for
$\za=\theta$, $\ze(\mu)=-1$ for $\za=\za_1\,, \za_2$ and the
justification of this choice will become clear below. 
The domain $E_p$ splits into two disjoint subsets $E_p^{(\pm)}$,
$E_p^{(+)} :=P_+^{p-2}\,,$ thus 
\eqn\sg{
\ze(\mu):=\pm 1\ \ {\rm for}\ \ \mu\in E_p^{(\pm)}\,.
}

Summarising  we are led to the following expression
for the quantised characters $\cc_{y}^{(p)}(\mu)\,, \, 
 y = \ol{y} \,t_{-\zl}\in \tC_p$:
\eqn\qfch{\eqalign{ 
\cc_{y}^{(p)}(\mu):
   &=
{e}^{{-2\pi\,i\, p\over 3} \zt(\mu)\, \zt(\zl)}\,
\Big( \bc_{\zl+\ol{y}^{-1}\cdot(0)}^{(p)}(\mu) 
     + \bc_{\zl+\ol{y}^{-1}\cdot(-2\fw_{1})}^{(p)}(\mu) 
     + \bc_{\zl+\ol{y}^{-1}\cdot(-2\fw_{2})}^{(p)}(\mu) \cr
   &+ (\FF(\mu)+2)
\big(\bc_{\zl+\ol{y}^{-1}\cdot(-\br)}^{(p)}(\mu)  
     + \bc_{\zl+\ol{y}^{-1}\cdot(-\fw_{2})}^{(p)}(\mu) 
     + \bc_{\zl+\ol{y}^{-1}\cdot(-\fw_{1})}^{(p)}(\mu)
\big)\Big)\,. 
}}
For $y\in \CC_p$ the overall  phase in \qfch\ disappears.  Taking
$\mu=0$ we define `$q$'-dimensions $D_y^{(p)}:=\cc_y^{(p)}(0)$
expressed by the `$q$'-dimensions of the integrable level $p-3$
case.

\PROP {\bf 6.1\ } {\it Let $x\,,\,y\in \CC_p\,,$ $\mu\in E_p\,.$ Then 
\eqn\qweyl{
  \cc_x^{(p)}(\mu)\ \cc_y^{(p)}(\mu)
  = \sum_{z\in \CC_p} {}^{(p)}N^{z}_{x,y} \ \cc_{z}^{(p)}(\mu)
\,, 
}
where
\eqn\qstrc{
{}^{(p)}N^{z}_{x,y} 
= \sum_{w'\in W^{[z\cdot k\zL_0]}} \,\det(\overline{w'})\ 
\ml^{x}_{w'\,z\, y^{-1}} 
= \sum_{w\in W^{[p]}} \, \det(\bw)\ \ml^{x}_{z w y^{-1}}
=  \sum_{w\in W^{[p]}} \, \det(\bw)\ N_{x, y }^{z w}
\,.
}
Furthermore the equality \qstrci\ holds true.}

\PROOF Since the  basic relations \fmult, \basr\ are conserved 
  the map $\cc_{y}\rightarrow \cc_{y}^{(p)}(\mu)$ is a ring
homomorphism,  so \weylI\ holds and it remains to use  \Vaa\ to
recover \qweyl, \qstrc. Finally the derivation of \qstrci\
parallels that of \strci\ using \afiii.
\endPROOF

The statement extends to $\cc_y^{(p)}(\mu)\,, \ y\in \tC_p$.
Given $ y\in \tC_p$ take $\zg^m\,y\in \CC_p$ with the appropriate
$m$. Then $\cc_y^{(p)}(\mu)= {e}^{{-2\pi\,ip\,m\,\tau(\mu)\over
3}}\,\cc_{\zg^m\,y}^{(p)}(\mu)$ and the product of characters
$\cc_y^{(p)}(\mu)\,, \ y\in \tC_p$ reduces to \qweyl, \qstrc\ due
to the symmetry ${}^{(3 p)} \bar{N}^{\iii(a\,
z)}_{\iii(a\,x)\,\iii(y)}= {}^{(3 p)}
\bar{N}^{\iii(z)}_{\iii(x)\,\iii(y)} \,, \, a\in A$, i.e., the
symmetries \invo, \invoi\  extend to ${}^{(p)}N^{z}_{x,y}$,
\eqn\qinvo{
   {}^{(p)}N_{x,y}^{z}={}^{(p)}N_{x,a y}^{a z}={}^{(p)}N_{x,z^*}^{y^*}
={}^{(p)}N_{x^*,y^*}^{z^*}\,, \quad \,, 
{}^{(p)}N_{x,y}^{\un}=\zd_{x,y^*}\,.
}
The action of the involution $*$
on the characters coincides with the complex conjugation
\eqn\chpr{
\cc_{y^{*}}^{(p)}(\mu)=\cc_y^{(p)}(\mu^{*})=
\cc_{y}^{(p)\,*}(\mu)\, (=\overline{\cc_y^{(p)}(\mu)})\,.
}
The second equality follows from $\ze(\mu)=\ze(\mu^*)$ and the
analogous equality for  the integrable characters.

 Using \afsc\ the first relation in \qinvo\ can be also 
rephrased in terms of
elements of $\CC_p$ only,  since $\cc_{_{\zg y}}^{(p)} =
\cc_{_{\zs_p(y)}}^{(p)}$ ($\zs_p(y)\,,$ $y\in \CC_p\,,  $ being the 
triality zero 
representative of $\zg y\in \tC_p\,$ on its $A_{[p]}$ orbit),
\eqn\afsca{
{}^{(p)}N_{x,\zs_p(y)}^{\zs_p(z)}
={}^{(p)}N_{x,y}^{z}\,, \qquad \cc_{_{\zs_p(\un)}}^{(p)}\, 
\cc_{_{y}}^{(p)} =
\cc_{_{\zs_p(y)}}^{(p)} 
\,.
}

The analogs of the basic examples in \fc\ read
\eqn\qfI{\eqalign{
\cc_{_{\zg}}^{(p)}(\mu)&=
{e}^{2 \pi i \, p \tau(\mu)\over 3}\,
(=\cc_{t_{-\zs_{[p-1]}(0)}}^{(p)}(\mu)=\cc_{\zg^2\,
t_{-\zs_{[p-1]}^2(0)}}^{(p)}(\mu)) =
\cc_{_{\zs_p(\un)}}^{(p)}(\mu)
\,,\cr
\cc_{_{w_0}}^{(p)}(\mu)&=2+ 
F(\mu)
= \bc_{_{(1,1)}}^{(3p)}(\mu)
 - R(\mu) - \bar{R}(\mu) + \ze(\mu)\,|R(\mu)|
\,,\cr
 \cc_{_{w_{20}}}^{(p)}(\mu)&=
{e}^{2 \pi i \, p \tau(\mu)\over
3}\, \cc_{_{t_{-\bL_1}}}^{(p)}(\mu)=
{e}^{2 \pi i \, p \tau(\mu)\over
3}\,\cc_{\zs_p^{-1}(w_{20})}^{(p)}(\mu)=
\bc_{_{(1,0)}}^{(p)}(\mu)
+1+ F(\mu)
\cr
&= \bc_{_{(3,0)}}^{(3p)}(\mu)
 - R(\mu) - \bar{R}(\mu) +\ze(\mu) \,|R(\mu)|
 \,.\cr
}}

In \qfI\ we have expressed  the characters in terms of the
integrable characters $\bc_{_{\iii(y)}}^{(3p)}(\mu)$ at (shifted)
level  $3p$. Since $E_p\subset P_+^{3p-3}$, taking $\mu\in E_p$
gives well defined expressions.  On the hyperplanes
$H_{\za}^{(p)}\cap E_p$ these characters reduce (up to a sign) to
the corresponding integrable characters $
\bc_{\iii(y)}^{(3 p)}(\mu)$ at (shifted) level  $3p$.
Indeed one proves

\LEMMA  {\bf 6.2\ } {\it Let $\mu\in E_p\cap \big(\cup_{\za\in
\bR_+}\,H_{\za}^{(p)}\big)$.
Then 
\eqn\rre{
r_{\ze}(\mu):= R(\mu) + \bar{R}(\mu) -\ze(\mu) \,|R(\mu)|=0\,
}
for $\ze(\mu)$ as in \sg.}
\PROOF
One easily checks  that for $\mu\in E_p\cap \big(\cup_{\za\in
\bR_+}\,H_{\za}^{(p)}\big)$ and $\ze(\mu)$ chosen as in \sg\
$R(\mu)$ can be cast into the form
$R(\mu)=-\ze(\mu)\,e^{{-2\pi \,i\over 3} \ze(\mu)}\, |R(\mu)|\,$
which implies the lemma. 
\endPROOF

The alternative expressions in \qfI\ representing the characters
$\cc_y^{(p)}(\mu)$ in terms of the integrable `$q$'-characters at
level $3p-3$ generalise to arbitrary $y\in \tC_p\,, \ \mu\in
E_p$. To simplify notation we shall omit the explicit dependence
on $\mu$ denoting the overall phase in \qfch\ by $\ze_y$.  Thus
for any $y=\by t_{\zl}\in \tC_p$ we obtain by straightforward
computation using  \froot, \rrela, \rrel, \qfI,  

\eqn\qfcha{
\cc_{y}^{(p)}
=\ze_y\,
\Big( \bc_{\iii(y)}^{(3 p)}  -r_\ze\,
\big(\bc_{\zl+\ol{y}^{-1}\cdot(-\br)}^{(p)}
     + \bc_{\zl+\ol{y}^{-1}\cdot(-\fw_{2})}^{(p)}
     + \bc_{\zl+\ol{y}^{-1}\cdot(-\fw_{1})}^{(p)}
\big)\Big)\,.
}
The second term in \qfcha\ admits
also a representation entirely in terms of integrable
`$q$'-characters $\bc_{\zn}^{(3 p)}(\mu) $ at level $3p-3$, with
weights $\zn\not \in$ Im$( \iii)$, using that
$\bc_{3\zl+2\br}^{(3p)}= r_{\ze}\, r_{-\ze}\,\bc_{\zl}^{(p)}$.
{}From Lemma 6.2 and the relations \rrela, \rrel\ it follows that
$r_{-\ze}(\mu)\not =0$ for any $\mu\in E_p$.
Finally we can also cast \qfcha\ into
 the form
\eqn\qcc{\eqalign{
\cc_{_{y}}^{(p)}(\mu)
=
\ze_y(\mu)\, \dd(\mu)
\Big(
\big(R^2 + F \bar{R}
\big)(\mu)& \sum_{w\in 
 \bar{A}}\, {e}^{{-2\pi i \over 3p}\,\la w(\iii(y)+\br) \,, \,
\mu+\br\ra}
\cr
-& \big( \bar{R}^2 +  F R\big)(\mu) \,
\sum_{w\in  \bar{A}}\, {e}^{{2\pi i\over 3p} \, \la
w(\iii(y)^*+\br)\,, \, \mu+\br\ra } \Big) \,.  
}} 

Lemma 6.2 and \qfcha\ imply 

\COR {\bf 6.3\ } {\it For any   $y\in  \CC_p$ and
 $\mu \in E_p\cap \big(\cup_{\za\in \bR_+} H_{\za}^{(p)}\big)\,,$}
\eqn\red{
\cc_{y}^{(p)}(\mu)=\bc_{\iii(y)}^{(3 p)}(\mu)\,. 
}

Despite of the relation \qstrci\ between  the structure constants
the product of characters $\cc_y^{(p)}(\mu)$ differs in general
from that of the integrable characters
$\bc_{\iii(y)}^{(3p)}(\mu)$ at level $3p-3$ since the
decomposition of the latter contains also terms
$\bc_{\zl}^{(3p)}(\mu)$ with $\zl\not \in$ Im$(\iii)$.  On the
other hand  the equality \qstrci\ together with \red\ -- the
latter property being enforced by the choice  \sg\ of the sign of
$F(\mu)$, require that on the intersection of the hyperplanes
$H_{\za}^{(p)}$ with $E_p$, the product   of the triality zero
integrable characters at shifted level $3p$ has to reduce to that
of the  characters \qfch.  Otherwise we run into contradiction,
i.e., the choice of sign \sg\ will appear to be inconsistent.
However it is easy to prove the above property of the standard
integrable characters at level $3p-3$, thus justifying a
posteriori the choice \sg. Namely we have

\LEMMA {\bf 6.4\ } {\it For  $ \mu\in P_+^{3p-3}\cap
\big(\cup_{\za\in \bR_+\,, \,  l\in \IZ} 
H_{\za}^{(lp)}\big)\,$ }
\eqn\infb{
\bc_{\iii(x)}^{(3p)}(\mu)\,\bc_{\iii(y)}^{(3p)}(\mu)
=\sum_{\zl\in {\rm Im}(\iii)}\,
{}^{(3p)}N_{\iii(x)\,\iii(y)}^{\zl}\  
\bc_{\zl}^{(3p)}(\mu)\,, \quad x,y\in \tC_p\,.
}

\PROOF 
The proof of the Lemma reduces to the proof of the following
property of the integrable characters at level $3p-3\,,$  $p>2$:

For $\mu\in P_+^{3p-3}\cap H_{\za}^{(lp)}\,,
\ \ \za\in \bR_+\,, \ l\in \IZ\,,$ and $\zl\in P_+^{3p-3}\,, \ 
\tau(\zl)=0\,,  \ \zl\not \in {\rm Im}(\iii)\,,$ 
\eqn\Vgg{
\bc_{\zl}^{(3p)}(\mu)=0\,.
}
If $\tau(\zl)=0\,,$ and $\zl\not \in $ Im$(\iii)$ then
$\zl+\br=3\,\zl'\,,$ for some 
$\zl'\in P_{+}^{p-3}+\br$. Hence 
$\bc_{\zl}^{(3p)}(\mu)={S_{\zl'-\br\,, \mu}^{(p)}\over 3\,
S_{0\,\mu}^{(3 p)}}$ and \Vgg\ follows from the vanishing of
$S_{\zl'-\br\,, \mu}^{(p)}$ for $\mu\in H_{\za}^{(lp)}\cap
P_+^{3p-3}\,.$
\endPROOF

\REMARK 
The case $p=2$ is degenerate (trivial) since the solutions of
\rre\ coincide with the whole $E_p$ and accordingly the triality
zero points in $P_+^{3p-3}\equiv P_+^{3}$ are all in Im($\iii$).
Hence the characters \qfch\ with $y=\un\,, \, w_{20}\,,
\, w_{10}\,, \, w_{0}\,$  coincide  with the corresponding
integrable characters at level $3p-3=3$ -- they realise the
triality zero fusion subalgebra at this level labelled by
$\{\zl=(0,0)\,, \,(3,0)\,, \, (0,3)\,, \, (1,1)\}\,.$ 
Thus the $\widehat{sl}(3)_k$  case $\kappa=k+3=3/2$ is analogous to
the $\widehat{sl}(2)_k$ case at  $\kappa=k+2=2/p\,, \ \  p$ --
odd, where the  admissible `$q$' - characters  $\cc_y^{(p)}(\mu)=
\bc_{\iii(y)}^{(2p)}(\mu)\,, \ y\in \CC_p\,,$
close the integer isospin ($\tau(\zl)=0$)   fusion subalgebra of the
$\widehat{sl}(2)$ integrable representations at shifted level $2p$;
the representative $ \CC_p$ of the admissible alcove is defined
as in \alcII, the latter formula being universal for any
$\widehat{sl}(n)_k$ and $k+n=n/p\,.$

%
%

\newsec{Pasquier--Verlinde type formula}

We have found $p^2$ vectors $\cc(\mu)=\{\cc_y^{(p)}(\mu)\,, \,
y\in \CC_p\}$ with $\mu\in E_p\,,$ which according to \qweyl\
provide eigenvectors common to all fusion matrices $N_y\,,
\ y\in \CC_p \,, $ $(N_y)_x^z={}^{(p)}N_{y, x}^z\,, $
and for any $y$ the numbers $\cc_y^{(p)}(\mu)\,\,$ are
eigenvalues of $N_y$ labelled by the set $E_p$.

\LEMMA {\bf 7.1\ } {\it Let $\mu\,, \mu'\in E_p\,.$ 
If $\cc(\mu)= \cc(\mu')$ then $\mu=\mu'$.}

\PROOF 
Recall that the  domain $E_p$ splits into two disjoint subsets
$E_p^{(\pm)}$ each being a subset of a fundamental domain in $P$
with respect to the shifted action of $W$ at level $p-3$.

{}From $\cc_{f_j}^{(p)}(\mu)= \cc_{f_j}^{(p)}(\mu')\,,\
j=0,1,2$  it follows that:

i) $\ze(\mu)=\ze(\mu')$,

\noindent
which  implies that both $\mu\,, \mu'
\in E_p^{(+)}$, or $\mu\,, \mu' \in E_p^{(-)}$, 

ii) $\bc_{\bL_i}^{(p)}(\mu)\,
=\bc_{\bL_i}^{(p)}(\mu')\,, \quad i=1,2\,,$ 

\noindent
which implies that $\mu'=\overline{w\cdot(\mu+(p-3)\zL_0)}\,, 
\, w\in W$. Hence  $\mu=\mu'$.
\endPROOF

Following standard arguments and taking into account the
properties \qinvo\ of the structure constants ${}^{(p)}N_{x,
y}^z$ the Lemma immediately leads to:
\COR  {\bf 7.2\ }
\eqn\funi{
\sum_{y\in \CC_p}\,\cc_y^{(p)}{(\mu)}\ \cc_{y}^{(p)\,*}{(\mu') }
=0\,, \quad \forall
\ \mu\,, \ \mu'\in E_p\,\ \  \mu\not =\mu'\,,
}
{\it and hence  $\{\cc(\mu)\,, \, \mu\in E_p\}$   is a linearly
independent set of (common) eigenvectors.}

Normalising the eigenvectors $\cc(\mu)$ (recall that 
$\cc_{\un}^{(p)}(\mu)=1$)
\eqn\rat{
\psi_y^{(\mu)}= \cc_y^{(p)}(\mu)\,\psi_{\un}^{(\mu)}\,, \quad
{1\over |\psi_{\un}^{(\mu)}|^2}=
\sum_{y\in \CC_p}\, |\cc_y^{(p)}(\mu)|^2\,, 
}
we can choose $\psi_{\un}^{(\mu)}$ real positive, so that
$\psi_y^{(\mu)\,*}=\psi_y^{(\mu^*)}=\psi_{y^*}^{(\mu)}\,,$ (see 
also \chpr{}).  Due to \funi\ the  square matrix $\psi_y^{(\mu)}$
is nonsingular and hence both its column and row vectors are
linearly 
independent. Thus we obtain a unitary matrix $\psi_y^{(\mu)}$,
\eqn\uni{
\sum_{y\in \CC_p}\,\psi_y^{(\mu)}\ \psi_y^{(\mu') *}=\zd_{\mu\,
\mu'}\,, \quad \sum_{\mu\in E_p}\,\psi_y^{(\mu)}\ 
\psi_x^{(\mu) *}=\zd_{y\, x}\,,
}
which diagonalises all $N_y$.  Indeed  using the second
(completeness) relation in \uni\ the formula  \qweyl\ converts
into \pasv, providing equivalent expression for the `$q$' analog
of the Weyl-Steinberg type formula \qstrc.  Hence  we recover the
Pasquier--Verlinde type formula for the fusion rule
multiplicities of the admissible representations at level
$k+3=3/p$ proposed in \FGP\ with now explicitly determined
eigenvector matrix $\psi_y^{(\mu)}$.

\medskip

A remaining technical problem is to 
perform explicitly the summation in
\rat. At least for  $\mu\in E_p\cap \big(\cup_{\za\in \bR_+}
H_{\za}^{(p)}\big)\,$ this  can be  easily done,  
 getting an explicit expression for the constant
$\psi_{\un}^{(\mu)}$, and hence for the corresponding matrix
elements of $\psi_y^{(\mu)}$, for this particular subset
of weights in $E_p$.  Indeed we have

\LEMMA  {\bf 7.3\ } {\it Let $y = \ol{y} \,t_{-\zl}\in \CC_p$ and
$\mu \in E_p\cap \big(\cup_{\za\in \bR_+} H_{\za}^{(p)}\big)\,.$
Then} 

\eqn\VIIa{\eqalign{
\psi_y^{(\mu)}&=
\sqrt{3}\, S_{\iii(y)\,\mu}^{(3p)}\,, \quad
\mu+\br\not = p\br\,; \cr
\psi_y^{(\mu)}&=\, S_{\iii(y)\,\mu}^{(3p)}\,, \quad \ \ \ 
\mu+\br = p\br\,.
}}
\PROOF According to \rat\ and \red\ it is sufficient to prove
the statement for $y=\un$. {}From \red, \Vgg,
it follows that for $\mu\in E_p\cap \big(\cup_{\za\in \bR_+}
H_{\za}^{(p)}\big)$ one has
$$
\sum_{y\in \CC_p}\, |\cc_y^{(p)}(\mu)|^2=\sum_{\zl\in
P_+^{3p-3}\,, \, \tau(\zl)=0}\, |\bc_{\zl}^{(3p)}(\mu)|^2
={1\over {(\sqrt{3}
S_{0\,\mu}^{(3p)})}^2}\, \sum_{l=0}^2\, \zd_{\mu,
\zs_{[3p-3]}^l(\mu)}.  
$$
The last  equality holds for any $\mu\in P_+^{3p-3}$ exploiting
standard properties of the modular matrices $S_{\zl\,\mu}^{(3p)}
$; see, e.g., \PZ.  Since the point $\mu+\br = p\br$ is a fixed
point for the $\zs_{[3p-3]}$ automorphism, the factor $\sqrt{3}$
does not appear in the second equality of \VIIa.
\endPROOF
According to the last Remark in the previous section in the case
$p=2$ the formul\ae\ \VIIa\ describe all matrix elements of the
eigenvector matrix $\psi_y^{(\mu)}$ and analogous formul\ae\
(with the factor $3$ substituted by $2$) hold for the whole
$sl(2)$ subseries at level $k+2=2/p\,.$

We conclude with the remark that the character ring constructed
here is an extension  of the ring of integrable  `$q$' -
characters at shifted level $p\,, $ with the two roots of
the quadratic polynomial \sqf{} of $F$.  The latter characters
are elements of the subring $\IZ[\omega]$ of the cyclotomic
extension $\IQ[\omega]$ of the rational numbers for
$\omega^{3p}=1$, see \BG.

\bigskip
\bigskip

{\bf Acknowledgements}
\medskip
We would like to thank V. Dobrev, V. Molotkov, Tch. Palev, I.
Penkov,  and J.-B. Zuber  for useful discussions, remarks, or,
suggestions. We also thank all colleagues who have shown interest
in this  work and/or have endured our explanations.  A.Ch.G.
thanks the A.v. Humboldt Foundation for financial support and the
Universities of Kaiserslautern and Bonn for hospitality.  V.B.P.
acknowledges the  financial support and hospitality of  INFN,
Sezione di Trieste, and the  partial support of the Bulgarian
National Research Foundation (contract $\Phi-643$).

\bigskip
{\bf Note added:}
\nref\GPW{Ganchev, A.Ch., Petkova, V.B.,  and Watts, G.M.T.,
to appear.}
\medskip

The multiplication rule encoded in the product of  the character
$\cc_{w_0}$  with any $\cc_{y}\,,$ see \fmt{}, has been
reproduced -- including the multiplicity two contribution, by an
explicit solution of the singular vectors decoupling equations at
generic level \GPW. Thus filling a gap in the computations in
\FGP\ all Pieri type formulae \fmt{}\ are now confirmed.

\bigskip
\listrefs
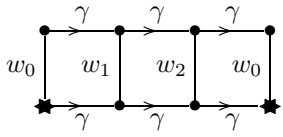
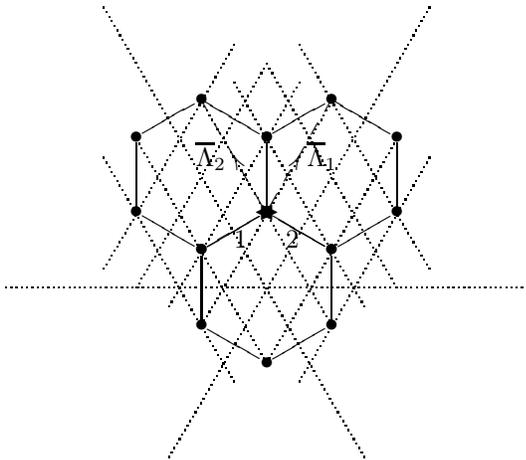
\bye